\documentclass{siamltex}
\usepackage[T1]{fontenc}
\usepackage{amsmath}%,amsthm}
\usepackage{amsfonts}
\usepackage{float}
\usepackage{graphicx}
\usepackage{colortbl}
\usepackage{subfigure}
\usepackage{color}

\definecolor{darkgreen}{rgb}{0.1,0.5,0.1}
\definecolor{darkblue}{rgb}{0.2,0.2,1.0}
\usepackage[colorlinks=true,linkcolor=darkblue,citecolor=darkblue,
            filecolor=darkblue,urlcolor=darkgreen]{hyperref}

\newenvironment{changebar}{}{}

\usepackage{url}
\usepackage[hmargin=1.0in]{geometry}

\begin{document}
\bibliographystyle{siam}

\newcommand{\mS}{\m{S}}
\newcommand{\mT}{\m{T}}
\newcommand{\ft}{\tilde{F}}
\newcommand{\qt}{\tilde{q}}
\newcommand{\tw}{\tilde{\omega}}

\newcommand{\TV}[1]{||#1||_\textup{TV}}

\newcommand{\qh}{\hat{q}}
\newcommand{\be}{\begin{equation}}
\newcommand{\ee}{\end{equation}}
\newcommand{\bq}{\mathbf{q}}
\newcommand{\bQ}{\mathbf{Q}}
\newcommand{\bx}{\mathbf{x}}
\newcommand{\br}{\mathbf{r}}
\newcommand{\imh}{{i-\frac{1}{2}}}
\newcommand{\Imh}{{I-\frac{1}{2}}}
\newcommand{\iph}{{i+\frac{1}{2}}}
\newcommand{\ipmh}{{i \pm \frac{1}{2}}}
\newcommand{\jmh}{{j-\frac{1}{2}}}
\newcommand{\jph}{{j+\frac{1}{2}}}
\newcommand{\Aop}{{\cal A}}
\newcommand{\Bop}{{\cal B}}
\newcommand{\Wop}{{\cal W}}
\newcommand{\Zop}{{\cal Z}}
\newcommand{\Wopt}{\widetilde{\cal W}}
\newcommand{\st}{\tilde{s}}
\newcommand{\Oop}{{\cal O}}
\newcommand{\DQ}{\Delta Q}
\newcommand{\Dq}{\Delta q}
\newcommand{\Dx}{\Delta x}
\newcommand{\Dy}{\Delta y}
\newcommand{\Du}{\Delta u}
\newcommand{\bu}{\mathbf{u}}
\newcommand{\bv}{\mathbf{v}}
\newcommand{\bw}{\mathbf{w}}
\newcommand{\bk}{\mathbf{k}}
\newcommand{\bU}{\mathbf{U}}
\newcommand{\bV}{\mathbf{V}}
\newcommand{\bF}{\mathbf{F}}
\newcommand{\Lop}{{\cal L}}
\newcommand{\Fop}{{\cal F}}
\newcommand{\Dofr}{{\cal D}(r)}
\newcommand{\Dt}{\Delta t}
\newcommand{\bbA}{\mathbb{A}}
\newcommand{\bbZ}{\mathbb{Z}}
\newcommand{\bbK}{\mathbb{K}}
\newcommand{\bM}{\mathbf{M}}
\newcommand{\bbI}{\mathbb{I}}
\newcommand{\bbb}{\mathbb{b}}
\newcommand{\bB}{\mathbf{B}}
\newcommand{\bR}{\mathbf{R}}
\newcommand{\bbe}{\mathbf{e}}
\newcommand{\bbone}{\mathbf{1}}

\newcommand{\amdQ}{\Aop^-\DQ}
\newcommand{\apdQ}{\Aop^+\DQ}
\newcommand{\apmdQ}{\Aop^\pm\DQ}
\newcommand{\amdq}{\Aop^-\Dq}
\newcommand{\apdq}{\Aop^+\Dq}
\newcommand{\amdqt}{\Aop^-\Delta {q}}
\newcommand{\apdqt}{\Aop^+\Delta {q}}
\newcommand{\apmdqt}{\Aop^\pm\Delta {q}}
\newcommand{\apmdq}{\Aop^\pm\Delta q}
\newcommand{\adqt}{\Aop\Delta q}
\newcommand{\adq}{\Aop\Delta q}
\newcommand{\dqt}{\Delta {q}}

\newcommand{\dx}{\Delta x}
\newcommand{\dt}{\Delta t}
\newcommand{\half}{\frac{1}{2}} 
\newcommand{\hfp}{\hat{f}_{j+\half}}
\newcommand{\hfn}{\hat{f}_{j-\half}}
\newcommand{\aik}{\alpha_{ij}}
\newcommand{\bik}{\beta_{ij}}
\newcommand{\lt}{\tilde{L}}
\newcommand{\ut}{\tilde{u}}

\newcommand{\hf}{\frac{1}{2}}
\newcommand{\fracStrut}{\rule[-1.0ex]{0pt}{3.1ex}}
\newcommand{\hfs}{\ensuremath{\frac{1}{2}}\fracStrut}
\newcommand{\scinot}[2]{\ensuremath{#1\times10^{#2}}}
\newcommand{\dee}{\mathrm{d}}
\newcommand{\dye}{\partial}
\newcommand{\diff}[2]{\frac{\dee #1}{\dee #2}}
\newcommand{\pdiff}[2]{\frac{\dye #1}{\dye #2}}
\newcommand{\Real}{\mathbb{R}}
\newcommand{\Complex}{\mathbb{C}}
% matrices
\newcommand{\m}[1]{\mathbf{#1}}
\newcommand{\mA}{\m{A}}
\newcommand{\mB}{\m{B}}
\newcommand{\mE}{\m{E}}
\newcommand{\mZ}{\m{Z}}
\newcommand{\mI}{\m{I}}
\newcommand{\mK}{\m{K}}
\newcommand{\mKt}{\m{\tilde{K}}}
\newcommand{\mL}{\m{L}}
\newcommand{\mQ}{\m{Q}}
\newcommand{\mP}{\m{P}}
\newcommand{\mPt}{\m{\tilde{P}}}
\newcommand{\mR}{\m{R}}
% these use the upgreek package to get non-italic greek, which doesn't
% seem to work with \mathbf so these have to be setup manually to
% match \m
% vectors
\renewcommand{\v}[1]{\boldsymbol{#1}}
\newcommand{\transpose}{^\mathrm{T}}
\newcommand{\bT}{\v{b}\transpose}
\newcommand{\vb}{\v{b}}
\newcommand{\vc}{\v{c}}
\newcommand{\vd}{\v{d}}
\newcommand{\ve}{\v{e}}
\newcommand{\vu}{\v{u}}
\newcommand{\vv}{\v{v}}
\newcommand{\vy}{\v{y}}
\newcommand{\vx}{\v{x}}
\newcommand{\vf}{\v{f}}
\newcommand{\vft}{\v{\tilde{f}}}
\newcommand{\vph}{\hat{\v{p}}}
\newcommand{\Matlab}{{\sc Matlab}\xspace}
\newcommand{\code}[1]{\textsf{#1}}
% the SSP coefficient
\newcommand{\sspcoef}{\mathcal{C}}
\newcommand{\bsspcoef}{\bf\mathcal{C}}

\newcommand{\clin}{\sspcoef_{\textup{lin}}}
\newcommand{\ceff}{\sspcoef_{\textup{eff}}}
\newcommand{\DtFE}{\Dt_{\textup{FE}}}
\newcommand{\Reff}{R_{\textup{eff}}}

\newcommand{\tL}{\tilde{\mL}}
\newcommand{\tF}{\tilde{F}}
\newcommand{\tR}{\tilde{R}}
\newcommand{\tr}{\tilde{r}}
\newcommand{\tz}{\tilde{z}}
\newcommand{\ta}{\tilde{a}}
\newcommand{\tk}{\tilde{k}}
\newcommand{\tp}{\tilde{p}}
\newcommand{\tb}{\tilde{b}}
\newcommand{\tc}{\tilde{\sspcoef}}
\newcommand{\tbeta}{\tilde{\beta}}
\newcommand{\tlambda}{\tilde{\lambda}}

\newcommand{\matalpha}{\boldsymbol{\upalpha}}
\newcommand{\matbeta}{\boldsymbol{\upbeta}}
\newcommand{\matmu}{\boldsymbol{\upmu}}
\newcommand{\matgamma}{\boldsymbol{\upgamma}}
\newcommand{\matlambda}{\boldsymbol{\uplambda}}
\newcommand{\mattheta}{\boldsymbol{\uptheta}}

\newcommand{\dd}[2]{\frac{\partial #1}{\partial #2}}

\newcommand{\Fig}[1]{Figure~\ref{fig:#1}}
\newcommand{\eqn}[1]{(\ref{#1})}
\newcommand{\Sec}[1]{Section~\ref{sect:#1}}
\newcommand{\Chap}[1]{Chapter~\ref{chap:#1}}

%Stegoton-specific macros
\newcommand{\rhomean}{\hat{\rho}}
\newcommand{\Kinvmean}{\hat{K}^{-1}}
\newcommand{\Kmean}{\hat{K}}
\newcommand{\Zmean}{\hat{Z}}
\newcommand{\cmean}{\hat{c}}

\newcommand{\shadeRow}{\rowcolor[rgb]{0.9, 0.9, 0.9}}
\newcommand{\shadeColumn}{\columncolor[rgb]{0.9, 0.9, 0.9}}

\title{High-Order Wave Propagation Algorithms for Hyperbolic Systems}
\author{
  David I. Ketcheson\thanks{King Abdullah University of Science and Technology,
    Box 4700, Thuwal, Saudi Arabia, 23955-6900.
   (\mbox{david.ketcheson@kaust.edu.sa}) } \and
   Matteo Parsani\thanks{King Abdullah University of Science and Technology,
    Box 4700, Thuwal, Saudi Arabia, 23955-6900.
   (\mbox{matteo.parsani@kaust.edu.sa}) }  \and 
   Randall J. LeVeque\thanks{Department of Applied Mathematics, University of Washington, Box 352420, Seattle, WA 98195-2420.
   (\mbox{rjl@uw.edu}) } }
\maketitle

\begin{abstract}
We present a finite volume method that is applicable to hyperbolic PDEs
including spatially varying and semilinear nonconservative systems. The spatial discretization, like
that of the well-known Clawpack software, is based on solving Riemann problems
and calculating fluctuations (not fluxes).  The implementation employs weighted
essentially non-oscillatory reconstruction in space and strong stability
preserving Runge-Kutta integration in time.  The method can be extended to
arbitrarily high order of accuracy and allows a well-balanced implementation
for capturing solutions of balance laws near steady state.  This
well-balancing is achieved through the $f$-wave Riemann solver and a novel
wave-slope WENO reconstruction procedure.  The wide
applicability and advantageous properties of the method are demonstrated
through numerical examples, including problems in nonconservative form,
problems with spatially varying fluxes, and problems involving near-equilibrium
solutions of balance laws.
\end{abstract}

\section{Introduction}
Many important physical phenomena are governed by hyperbolic systems of 
conservation laws. In one space dimension the standard conservation law has the form
\be
\label{conslaw}
q_t + f(q)_x = 0,
\ee
where the components of $q \in \mathbb{R}^m$ are conserved quantities and 
the components of $f:\mathbb{R}^m \times \mathbb{R}^m \rightarrow \mathbb{R}^m$  are
the corresponding fluxes. 
Very many numerical methods have been developed for the solution of 
\eqref{conslaw}; some of the most successful are the high resolution
Godunov-type methods based on the use of Riemann solvers and nonlinear limiters.
These and other methods are generally based on {\em flux-differencing} and
make explicit use of the flux function $f$.

Herein we also consider systems with spatially varying flux:
\be
\label{varflux}
\kappa(x)q_t + f(q,x)_x  = \psi(q,x), 
\ee
and spatially varying linear systems not in conservation form:
\be
\label{varlin}
\kappa(x)q_t + A(x)q_x  = \psi(q,x), 
\ee
as well as their two-dimensional extensions.
%Examples of (\ref{conslaw}-\ref{hypsys}) include acoustics and elasticity 
%in heterogeneous media. 
Wave-propagation methods of up to second-order accuracy have been developed 
for such systems in, e.g. \cite{leveque1997,leveque2002}.  These methods are
based on wave-propagation Riemann solvers, which compute {\em fluctuations}, 
(i.e. traveling discontinuities) rather than fluxes and thus can be applied 
to \eqref{varflux} or \eqref{varlin}
just as easily as to the conservation law \eqref{conslaw}.  
%Like those methods,
%the methods in this work could in principle be applied to nonlinear, nonconservative
%systems; see, e.g. \cite{castro2008b} for a discussion of the difficulties
%arising in such systems.

Second-order methods may often be the best choice in terms of a balance 
between computational cost and desired resolution, especially for problems 
with solutions dominated by shocks or other discontinuities with relatively 
simple structures between these discontinuities. For problems containing 
complicated smooth solution structures, where the accurate 
resolution of small scales is required (e.g. simulation of compressible 
turbulence, computational aeroacoustics, computational 
electromagnetism, turbulent combustion etc.), schemes with higher 
order accuracy are desirable.

The purpose of this work is to present a numerical method that combines the
advantages of wave-propagation solvers with high order of accuracy.
The basic discretization approach was presented already in \cite{ketcheson2006};
here, we give a more
detailed presentation and demonstrate the wide applicability of the method.
The new method combines the notions of wave propagation 
(\cite{leveque1997,levequefvbook}) and
the method of lines, and can in principle be extended to arbitrarily
high order accuracy by the use of high order accurate spatial reconstructions
and high order accurate ordinary differential equation (ODE) solvers. 
The implementation presented here is based on the fifth-order accurate weighted essentially non-oscillatory (WENO) reconstruction and
a fourth-order accurate strong-stability-preserving Runge-Kutta (RK) scheme.
%The wave propagation is particularly useful in extending the great variety of existing spatial discretization techniques to hyperbolic system that are not in conservation form, e.g. system \ref{hypsys}. 
We restrict our attention
to problems in one or two dimensions, although the method may be extended
in a straightforward manner to higher dimensions.

Although the method described here can be applied to classical hyperbolic
systems \eqref{conslaw}, in that case it is equivalent
to a standard finite volume WENO flux-differencing scheme,
as long as component-wise reconstruction and a conservative wave propagation
Riemann solver (such as a Roe or HLL solver) are used.  
For this reason, we focus on problems in non-conservative form or with explicit
spatial dependence in the flux or source terms, which can be challenging
for traditional discretizations.

Another approach to high order discretization of hyperbolic
PDEs, referred to as the ADER method, has
been developed by Titarev \& Toro \cite{titarev2002} and
subsequent authors.  That approach uses
the Cauchy-Kovalewski procedure and has the advantage of leading to
one-step time discretization.  The method of lines approach used in the 
present work seems more straightforward and
allows manipulation of the method's properties by the use of different
time integrators, but requires the evaluation of multiple stages per
time step.

A similar class of conservative, well-balanced, and high-order accurate methods
%, applicable to nonconservative hyperbolic systems, 
has been developed by Castro, Gallardo, Pares,
and their coauthors; see, e.g. \cite{castro2008}.  Those methods also use 
WENO reconstruction and Runge-Kutta time stepping in conjunction with Riemann
solvers, and lead to a discretization with a form similar to that presented here and in
\cite{ketcheson2006}.  Those methods
have recently been combined with the ADER approach; see \cite{dumbser2009}.
The present method differs in the approach to reconstruction and the kind of
Riemann solvers used.  These differences result in some useful features:
our method also handles systems \eqref{varflux}-\eqref{varlin} with capacity function $\kappa$ 
and it can make use of $f$-wave Riemann solvers \cite{bale2002}
as well as wave-slope reconstruction and achieve high order convergence even
for some problems with discontinuous coefficients.  Finally, the method can
immediately be applied to very many interesting problems because the implementation
is based on Clawpack Riemann solvers, which are available for a great variety
of hyperbolic systems.  A potential drawback of the current implementation of our method
is that, for well-balancing, it relies on discretizing the source term such that
its effect is collocated at cell interfaces only.

The methods described in this paper are implemented in the software package
SharpClaw, which is freely available on the web at
\url{http://www.clawpack.org}.  SharpClaw employs the same interface that is
used in Clawpack \cite{levequefvbook} for problem specification and setup, as
well as for the necessary Riemann solvers.  This makes it simple to apply
SharpClaw to a problem that has been set up in Clawpack.  These methods
have also been incorporated into PyClaw\cite{PyClaw-2011-SISC}, which allows them to be run in
parallel on large supercomputers.

The paper is organized as follows.  In \Sec{god-waveprop}, we present Godunov's method for
linear hyperbolic PDEs in wave propagation form \cite{levequefvbook}.  This method is 
extended to high order in \Sec{ho-waveprop} by introducing a high-order
reconstruction based on cell averages. Generalization to spatially varying 
nonconservative linear systems is presented in \ref{sect:svnonconshypsys}. 
Further extensions and details of the method are presented in the remainder of 
Section 2.  Numerical examples, including application to acoustics, elasticity, 
and shallow water waves,are presented in Section \ref{sect:applications}.

\section{Semi-discrete wave-propagation}
The wave-propagation algorithm was first introduced by LeVeque \cite{leveque1997} in 1997 in the framework of high resolution finite volume methods for solving hyperbolic systems of equations. The scheme is conservative, second-order accurate in smooth regions, and captures
shocks without introducing spurious oscillations.  In this section, we extend
the wave-propagation algorithm to high order accuracy through
use of high order reconstruction and time marching.
For simplicity, we focus on the one dimensional scheme and then briefly describe the extension to
two dimensions.

\subsection{Riemann Problems and Notation}
The notation for Riemann solutions used in this paper comes primarily 
from \cite{levequefvbook}, and is motivated by consideration of the
linear hyperbolic system
\be \label{eq:linhypsys}
q_t + A q_x = 0.
\ee
Here $q\in\mathbb{R}^m$ and $A\in\mathbb{R}^{m\times m}$.  System \eqref{eq:linhypsys}
is said to be hyperbolic if $A$ is diagonalizable with real
eigenvalues; we will henceforth assume this to be the case.
Let $s^p$ and $r^p$  for $1\le p\le m$ denote the eigenvalues and the corresponding 
right eigenvectors of $A$ with the eigenvalues ordered so that 
$s^1 \leq s^2 \leq \ldots \leq s^m$. 

\begin{figure}
\centering
\includegraphics[width=4in]{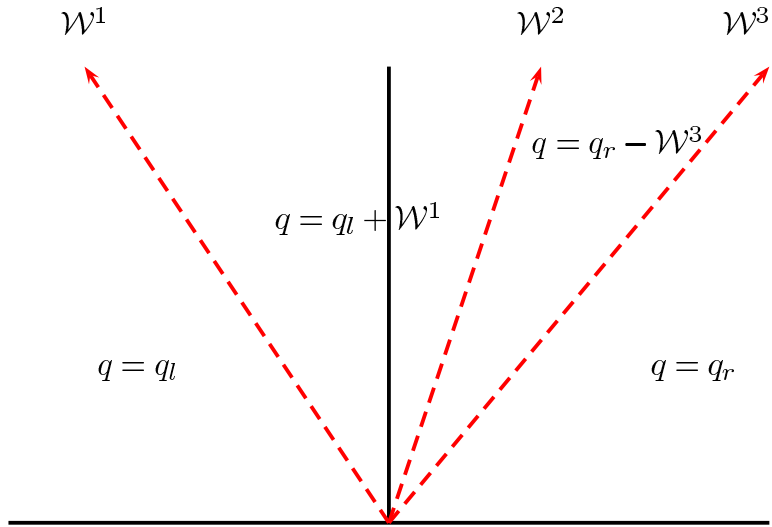}
%\begin{pspicture}(0,0)(10,7)
%\psset{unit=0.85cm}
%\psset{linecolor=black,linewidth=0.05}
%\psline(0,0)(10,0)
%\psline(5,0)(5,6)
%\psline[linecolor=red,linestyle=dashed,arrows=->](5,0)(1,6)
%\rput[c](1,6.6){$\Wop^1$}
%\psline[linecolor=red,linestyle=dashed,arrows=->](5,0)(7,6)
%\rput[c](7,6.6){$\Wop^2$}
%\psline[linecolor=red,linestyle=dashed,arrows=->](5,0)(10,6)
%\rput[c](9.7,6.6){$\Wop^3$}
%\rput[c](1.7,2.0){$q=q_l$}
%\rput[c](9,2.0){$q=q_r$}
%\rput[c](4.6,4.0){$q=q_l+\Wop^1$}
%\rput[c](8,5){$q=q_r-\Wop^3$}
%\end{pspicture}
\caption{The wave propagation solution of the Riemann problem.  The horizontal
axis corresponds to space and the vertical axis to time.\label{fig:waveprop}}
\end{figure}

Consider the {\em Riemann problem} consisting of
\eqref{eq:linhypsys} together with initial data
\begin{align} \label{eq:riemann-data}
q(x,0)=\left\{ \begin{array}{cc} q_l & x<0 \\ q_r & x>0 \end{array}\right. .
\end{align}
The solution for $t>0$ is piecewise constant with $m$ discontinuities, the $p$th one proportional
to $r^p$ and moving at speed $s^p$.  They can be obtained by decomposing the 
difference $q_r-q_l$ in terms of the eigenvectors $r^p$:
\begin{align} \label{eq:wavedef}
q_r-q_l=\sum_p \alpha^p r^p = \sum_p \Wop^p.
\end{align}
We refer to the vectors $\Wop^p$ as waves.  Each wave is a jump
discontinuity along the ray $x=s^p t$.  An example
solution is pictured in Figure \ref{fig:waveprop} for $m=3$.
For brevity, we will sometimes refer to the Riemann problem with initial left state
$q_l$ and initial right state $q_r$ as {\em the Riemann problem with initial states 
$(q_l,q_r)$.}

In a finite volume method, it is useful to define notation for the net effect of
all left- or right-going waves:
\begin{subequations}\label{eq:fluctdef}
\begin{eqnarray} 
\amdq  & \equiv & \sum_{p=1}^m \left(s^p\right)^- \Wop^p\\
\apdq  & \equiv & \sum_{p=1}^m \left(s^p\right)^+ \Wop^p.
\end{eqnarray}
\end{subequations}
Here and throughout, $(x)^\pm$ denotes the positive or negative part of $x$:
$$(x)^-=\min(x,0) \quad \quad (x)^+=\max(x,0).$$
The symbols $\apmdq$, referred to as {\em fluctuations},
should be interpreted as 
single entities that 
represent the net effect of all waves travelling to the right or left.
% RJL added:
The notation is motivated by the constant coefficient linear system
\eqref{eq:linhypsys},  in which case $\apmdq = A^{\pm}(q_r-q_l)$, where $A^-$
(respectively $A^+$) is the matrix obtained by setting all postive (respectively
negative) eigenvalues of $A$ to zero. See \cite{leveque1997} or
\cite{levequefvbook} for more details.

The notation for waves and fluctuations defined in \eqref{eq:wavedef}
and \eqref{eq:fluctdef} can also be used to describe numerical solutions of
Riemann problems for nonlinear systems if the numerical solver approximates
the solution by a series of propagating jump discontinuities, which is
very often the case.  Because the approximate Riemann solution for a nonlinear
system depends not only on the difference $q_r-q_l$ but on the values
of the states, we will sometimes employ for clarity the notation 
$\Wop^p(q_l,q_r)$ to denote the
$p$th wave in the solution of the Riemann problem with initial states $(q_l,q_r)$.
% RJL added more:
\begin{changebar}
In this case the vectors $r^p$ used for the decomposition
\eqn{eq:wavedef} are typically eigenvectors of an averaged Jacobian
matrix $\bar A$, and the $s^p$ are the corresponding wave speeds.
Then $\apmdq = \bar A^{\pm}(q_r-q_l)$.  But there are other possible
ways to define the $r^p$, $s^p$ and hence the fluctuations.  For
example, using the HLL approximate Riemann solver \cite{ha-la-vl}
would always use only two waves regardless of the size of the system.
Or, for the case of spatially varying flux functions, the left-going
waves may be defined by eigenvectors of the Jacobian $f'(q_l)$ while
the right-going waves are defined by eigenvectors of the Jacobian
$f'(q_r)$. One could combine these to create a matrix $\bar A$
having this set of eigenvectors and the corresponding eigenvalues,
but this is not required for implementation of the method.  For
these reasons we use the more general notation $\apmdq$ for the
fluctuations defined by \eqref{eq:fluctdef}.
\end{changebar}

%=====================================================================
\subsection{First-order Godunov's method\label{sect:god-waveprop}}
%=====================================================================
%Consider the {\em Riemann problem} consisting of
%\eqref{eq:linhypsys} together with initial data
%\begin{align*} 
%q(x,0)=\left\{ \begin{array}{cc} q_l & x<0 \\ q_r & x>0 \end{array}\right.
%\end{align*}
%The solution for $t>0$ consists of $m$ discontinuities, the $p$th one proportional
%to $r^p$ and moving at speed $s^p$.  They can be obtained by decomposing the 
%difference $q_r-q_l$ in terms of the eigenvectors $r^p$:
%\be
%\label{eq:waveDecompTh}
%q_r-q_l=\sum_p \alpha^p r^p = \sum_p \Wop^p.
%\ee
%The coefficients $\alpha^p$ are given by
%\be
%\label{eq:alphaCoeff}
%\alpha^p = l^p(q_r-q_l),
%\ee
%where $l^p$ are the left eigenvectors are the matrix $A$ \cite{levequefvbook}. We refer to the vectors $\Wop^p$ as waves.  Each wave is a jump
%discontinuity along the ray $x=s^p t$ in phase space.  The solution is pictured in Figure \ref{fig:waveprop} for $m=3$.

Consider the constant-coefficient linear system in one dimension 
\eqref{eq:linhypsys}.
Taking a finite volume approach, we define the cell averages (i.e. the solution variables)
\begin{align*}
Q_i(t) = \frac{1}{\Dx}\int_{x_\imh}^{x_\iph} q(x,t) dx,
\end{align*}
where the index $i$ and the quantity $\Dx$ denote the cell index and the cell size, respectively.
To solve the linear system \eqref{eq:linhypsys}, we initially approximate the solution
$q(x,t)$ by these cell averages; that is, at $t=t_0$ we define the piecewise-constant
function
\be \label{eq:pwc}
\qt(x,t_0) = Q_i \mbox{ for } x\in (x_\imh,x_\iph).
\ee
The linear system \eqref{eq:linhypsys} with initial data $\qt$ consists locally of a series 
of Riemann problems, one at each interface $x_\imh$.  
%Therefore, the solution variables are updated in each time step using the wave structure determined by solving Riemann problems at cell edges.  
The Riemann problem at $x_\imh$ consists of (\ref{eq:linhypsys}) with the piecewise constant initial data
\begin{align*} 
q(x,0)=\left\{\begin{array}{ll} Q_{i-1} & x<x_\imh \\ 
                                                Q_i & x>x_\imh 
                     \end{array}\right. .
\end{align*}
%Let $s^p$ and $r^p$  for $1\le p\le m$ denote the eigenvalues and the corresponding 
%right eigenvectors of $A$ with the eigenvalues ordered so that 
%$s^1 \leq s^2 \leq \ldots \leq s^m$. 
As discussed above, the solution of the Riemann problem is expressed as a set of waves obtained by decomposing
the jump in $Q$ in terms of the eigenvectors of $A$:
\be
\label{eq:waveDecompFV}
Q_i-Q_{i-1}=\sum_p \alpha^p_\imh r^p_\imh = \sum_p \Wop^p(Q_{i-1},Q_i).
\ee

Let $\qt(x,t_0+\Dt)$ denote the exact evolution of $\qt$ after a time
increment $\Dt$.  If we take $\Dt$ small enough that the waves from
adjacent interfaces do not pass through more than one cell, then 
we can integrate \eqref{eq:linhypsys} over $[x_\imh,x_{\iph}] \times [0,\Dt]$
and divide by $\Dx$ to obtain 
\begin{align} \label{eq:1int}
Q_i(t_0+\Dt) - Q_i(t_0) & = -\frac{1}{\Dx} \int_{x_\imh}^{x_\iph} A \, \qt_x(x,t_0+\Dt) dx.
\end{align}
Here $\qt_x$ should be understood in the sense of distributions.

\begin{figure}
\centering
\includegraphics[width=4in]{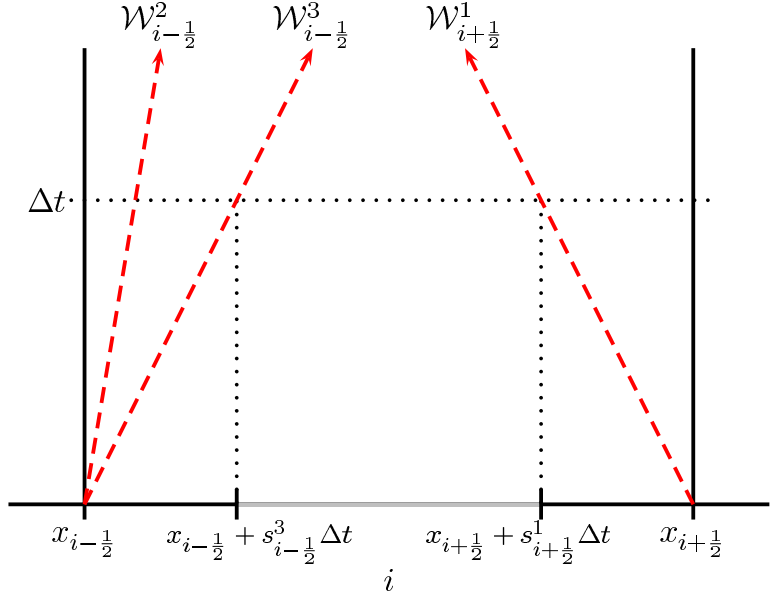}
%\begin{pspicture}(0,-1)(10,7)
%\psset{unit=0.85cm}
%\psset{linecolor=black,linewidth=0.05}
%\psline(0,0)(10,0)
%\psline[linecolor=lightgray](3,0)(7,0)
%\psline(1,-0.2)(1,6)
%\psline(9,-0.2)(9,6)
%\rput[c](1,-0.5){$x_\imh$}
%\rput[c](9,-0.5){$x_\iph$}
%\psline[linestyle=dotted](0.8,4)(9.2,4)
%\rput[c](0.5,4){$\Dt$}
%\psline[linestyle=dotted](3,0)(3,4)
%\psline(3,-0.2)(3,0.2)
%\rput[c](3.3,-0.5){\footnotesize $x_\imh+s^3_\imh\Dt$}
%\rput[c](5.0,-1.0){\normalsize $i$}
%\psline[linestyle=dotted](7,0)(7,4)
%\psline(7,-0.2)(7,0.2)
%\rput[c](6.7,-0.5){\footnotesize $x_\iph+s^1_\iph\Dt$}
%\psline[linecolor=red,linestyle=dashed,arrows=->](9,0)(6,6)
%\rput[c](6,6.3){$\Wop^1_\iph$}
%\psline[linecolor=red,linestyle=dashed,arrows=->](1,0)(2,6)
%\rput[c](2,6.3){$\Wop^2_\imh$}
%\psline[linecolor=red,linestyle=dashed,arrows=->](1,0)(4,6)
%\rput[c](4,6.3){$\Wop^3_\imh$}
%\end{pspicture}
\caption{Time evolution of the reconstructed solution $\qt$ in cell $i$.\label{fig:qt-evol}}
\end{figure}

We can split the integral above into three parts, representing the
Riemann fans from the two interfaces, and the remaining piece:
\begin{align} \label{eq:3ints}
\int_{x_\imh}^{x_\iph} A \qt_x dx & = \int_{x_\imh}^{x_\imh+s^R \Dt} A \qt_x dx + \int_{x_\iph+s^L \Dt}^{x_\iph} A \qt_x dx + \int_{x_\imh+s^R \Dt}^{x_\iph+s^L \Dt} A \qt_x dx.
\end{align}
The relevant regions are depicted in Figure \ref{fig:qt-evol}.  Here we have defined
$s^L=\min(s^1_\iph,0)$ and $s^R=\max(s^m_\imh,0)$.
The third integral in \eqref{eq:3ints} vanishes because $\qt(x,\Dt)$ is 
constant outside the Riemann fans, by the definition \eqref{eq:pwc}.  Hence
\eqref{eq:3ints} reduces to
\begin{subequations}
\label{eq:3intsb}
\begin{align} 
\int_{x_\imh}^{x_\iph} A \qt_x dx & = \Dt \sum_{p=1}^m \left(s^p_\imh\right)^+ \Wop^p_\imh \label{eq:3ints-waveDec}
    +\Dt \sum_{p=1}^m \left(s^p_\iph\right)^- \Wop^p_\iph \\
 & = \Dt\left(\Aop^+\DQ_\imh+\Aop^-\DQ_\iph\right),\label{eq:3ints-fluct}
\end{align}
\end{subequations}
%where $(x)^\pm$ denotes the positive or negative part of $x$:
%$$(x)^-=\min(x,0),  \quad \quad (x)^+=\max(x,0),$$  
where the {\em fluctuations} $\Aop^-\DQ_\iph$ and $\Aop^+\DQ_\imh$ are defined by
\begin{subequations}\label{pfluct}
\begin{eqnarray} 
\Aop^-\DQ_\iph  & \equiv & \sum_{p=1}^m \left(s^p_\iph\right)^- \Wop^p_\iph, \qquad \Wop^p_\iph \equiv \Wop^p(Q_{i},Q_{i+1})\\
\Aop^+\DQ_\imh  & \equiv & \sum_{p=1}^m \left(s^p_\imh\right)^+ \Wop^p_\imh, \qquad \Wop^p_\imh \equiv \Wop^p(Q_{i-1},Q_i).
\label{mfluct} \end{eqnarray}
\end{subequations}
Note again that the fluctuations $\Aop^+\DQ_\imh$ and $\Aop^-\DQ_\iph$
are motivated by the idea of a matrix-vector product but 
should be interpreted as single entities that 
represent the net effect of all waves travelling to the right or left.
The upper-case $Q$ in the fluctuations is meant to emphasize that they
are based on differences of cell averages.
For instance, the fluctuation $\Aop^+\DQ_\imh$ corresponds to the effect
of right-going waves from the Riemann problem with initial states $(Q_{i-1},Q_i)$.
%This is illustrated in \Fig{qt-evol} for $m=3$.  
%The waves $\Wop_\imh$ and speeds $s_\imh$ are those resulting
%from the Riemann problem with left and right states $Q_{i-1}$ and
%$Q_i$, respectively. 

Substituting \eqref{eq:3ints-fluct} into \eqref{eq:1int}, we obtain the scheme
$$
Q_i^{n+1}-Q_i^n = -\frac{\Dt}{\Dx} 
\left(\Aop^+\DQ_\imh + \Aop^-\DQ_\iph \right).
$$
Dividing by $\Dt$ and taking the limit as $\Dt$ approaches zero,
we obtain the semi-discrete wave-propagation form of the (first-order) Godunov's scheme
\be
\label{sdform_1}
\frac{\partial Q_i}{\partial t} = -\frac{1}{\Dx} 
\left(\Aop^+\DQ_\imh + \Aop^-\DQ_\iph \right).
\ee
Equation \eqref{sdform_1} constitutes a linear system of ordinary differential equations (ODEs) that may be 
integrated, for instance, with a Runge-Kutta method. It is clear from the derivation that this
scheme reduces to the corresponding flux-differencing scheme when applied
to systems written in conservation form, e.g. system (\ref{varflux}). 
%The advantage of the proposed scheme over flux-differencing schemes lies in the
%ability to solve systems which are not in conservation form, e.g. the general quasilinear variable-coefficient hyperbolic system (\ref{hypsys}). Since systems of this form generally cannot be rewritten in terms of a flux function, fluctuations are calculated in terms of the decomposition (\ref{eq:3ints}).

%=====================================================================
\subsection{Extension to higher order\label{sect:ho-waveprop}}
%=====================================================================
The method of the previous section is only first-order accurate
in space.  In order to improve the spatial accuracy, we replace
the piecewise-constant approximation \eqref{eq:pwc} by a 
piecewise-polynomial approximation that is accurate to order
$p$ in regions where the solution is smooth:
\begin{align} \label{eq:pwp}
\qt(x,t_0) & = \qt_i(x) \mbox{ for } x\in (x_\imh,x_\iph),
\end{align}
where
$$
\qt_i(x)=q(x,t_0)+\Oop(\Dx^{p+1}).
$$
%The reconstruction is illustrated in \Fig{piecewise-poly}.

%\begin{figure}
%\centering
%    \begin{pspicture}(0,0.5)(8,5)
%      \psset{unit=0.85cm}
%      \psset{linewidth=0.05}
%      \psline(0,1.5)(8,1.5) 
%      \psdots(0,1.5)(2,1.5)(4,1.5)(6,1.5)(8,1.5)
%      \rput[c](4,1){$x_\imh$}
%      \psline[linecolor=red,linestyle=dashed](0,2)(2,2)(2,3)(4,3)(4,4.4)(6,4.4)(6,3.8)(8,3.8)
%      \pscurve(2,2.4)(3,2.8)(4,3.8)
%      \pscurve(4,4.2)(5,4.6)(6,3.8)
%    \end{pspicture}
%\caption{Illustration of piecewise polynomial reconstruction from cell averages.\label{fig:piecewise-poly}}
%\end{figure}

Integration of $A\qt_x$ over $[x_\imh,x_{\iph}]$
again yields \eqref{eq:3ints}, but now the third integral is non-zero in general,
%divide by $\Dt\Dx$, and take the limit as $\Dt$ approaches zero.
%We now find that the third integral in \eqref{eq:3ints} contributes,
since $\qt$ is not constant outside the Riemann fans.  Define
\begin{changebar}
\begin{align}
q^R_\imh \equiv \lim_{x \to x_\imh^+} \qt_i(x) \qquad \qquad q^L_\iph \equiv \lim_{x \to x_\iph^-} \qt_i(x),
\end{align}
\end{changebar}
where superscripts $L$ and $R$ refer respectively to the left and the right state of the interface considered.  Then in place of \eqref{eq:3intsb}, we now obtain (neglecting terms of order $\Oop(\Dt^2)$)
\begin{subequations}
\label{eq:3intsc}
\begin{align} 
\int_{x_\imh}^{x_\iph} A \qt_x dx & \approx \Dt \sum_{p=1}^m \left(s^p_\imh\right)^+ \Wop^p_\imh 
    +\Dt \sum_{p=1}^m \left(s^p_\iph\right)^- \Wop^p_\iph+ \int_{x_\imh+s^R\Dt}^{x_\iph+s^L\Dt} A \qt_x dx \label{eq:3intsc-waveDec}\\
 & = \Dt\left(\apdqt_\imh+\amdqt_\iph\right) + A(q^L_\iph-q^R_\imh).\label{eq:3intsc-fluct}
\end{align}
\end{subequations}

The resulting fully-discrete scheme is thus
$$
Q_i^{n+1}-Q_i^n = -\frac{\Dt}{\Dx} 
\left(\apdqt_\imh+ \amdqt_\iph + A(q^L_\iph-q^R_\imh)\right).
$$
We use the notation $\apmdqt$ instead of $\Aop^\pm\DQ$ because the states
in the Riemann problems are not the cell averages, but rather the reconstructed
interface values.  
%For the same reason, we use $\st,\Wopt$ to denote the speeds
%and waves resulting from these Riemann problems.  
In other words, the
fluctuations at $x_\imh$ are defined by
\begin{align*}
\apmdqt_\imh & = \sum_{p=1}^m \left(s^p(q^L_\imh,q^R_\imh)\right)^\pm \Wop^p(q^L_\imh,q^R_\imh).
\end{align*}
For instance, the fluctuation $\apdqt_\imh$ corresponds
to the effect of right-going waves from the Riemann problem
with initial states $(q^L_\imh,q^R_\imh)$. Moreover, we can view
the term $A(q^L_\iph-q^R_\imh)$ as the sum of both the left-
and right-going fluctuations resulting from a Riemann problem with initial states
$(q^R_\imh,q^L_\iph)$.
It is natural to denote this term,
which we refer to as a {\em total fluctuation}, by $\adqt_i$:
\begin{align*}
\adqt_i & = \sum_{p=1}^m \left(s^p(q^R_\imh,q^L_\iph)\right)^\pm \Wop^p(q^R_\imh,q^L_\iph).
\end{align*}
Dividing by $\Dt$ and taking the limit as $\Dt$ approaches zero,
we obtain the semi-discrete scheme
\be
\label{sdform_1a}
\frac{\partial Q_i}{\partial t} = -\frac{1}{\Dx} 
\left(\apdqt_\imh + \amdqt_\iph + \adqt_i\right).
\ee

%%=====================================================================
%\subsection{Variable Coefficient Linear Systems}
%    \label{sect:varcoeff}
%%=====================================================================
%We now generalize the method to solve linear hyperbolic systems with 
%variable coefficients:
%\be \label{eq:linhypsys-varcoeff}
%q_t + A(x) q_x = 0.
%\ee
%We assume that the system is hyperbolic for all $x$ and that
%$A(x)$ is piecewise-constant, with points of
%discontinuity aligned with grid interfaces.  Thus $A(x)$ is
%given by a constant matrix $A_i$ within grid cell $i$.  The Riemann
%problem at $x_\imh$ is now given by
%\eqref{eq:linhypsys-varcoeff} together with 
%\be q(x,0)=\left\{ \begin{array}{cc} q^-_\imh & x<x_\imh \\ q^+_\imh & x>x_\imh \end{array}\right.
%\quad \quad
%A(x)=\left\{ \begin{array}{cc} A_{i-1} & x<x_\imh \\ A_{i} & x>x_\imh \end{array}\right.
%\ee
%As in the constant coefficient case, the Riemann solution consists
%of waves, but now the left-going waves are multiples of the eigenvectors
%$r_{i-1}$ of $A_{i-1}$ while the right-going waves are multiples
%of the eigenvectors $r_i$ of $A_i$.  Thus the semi-discrete
%scheme is again \eqref{sdform_1a}, but with fluctuations corresponding
%to these waves and with $A=A_i$.

%=====================================================================
\subsection{Spatially varying linear systems\label{sect:svnonconshypsys}}
%=====================================================================
Next we generalize the method to solve one-dimensional spatially varying 
nonconservative linear systems:
\be \label{eq:nlhypsys}
q_t + A(x) q_x = 0.
\ee
We again assume that $A$ is a constant function of $x$ within
each cell, so we can write $A(x)=A_i(q)$.  In the special case
that $A$ is the Jacobian matrix of some function $f$, 
\eqref{eq:nlhypsys} corresponds to the conservation law \eqref{conslaw}.
Our method can be applied to the general system \eqref{eq:nlhypsys}
as long as the physically meaningful solution to the Riemann problem can be approximated.
\begin{changebar}
This may be nontrivial for a nonconservative nonlinear problem, as
discussed in many papers such as
\cite{ca-lef:sw,ca-di-ch-etal:wb,co-le-no-pe:multdist,dm-lef-mu,lef-tza:weak}.
We do not go into this here, but assume that our numerical approach is to be
applied to a problem for which the user has a means of solving the Riemann
problem in terms of 
discontinuities (waves) $\Wop^p_\imh$ propagating at speeds $s^p_\imh$, and
hence can define fluctuations. This is
the case for any linearized Riemann solver, and often for other approximate solvers.
\end{changebar}
Then the scheme is given by
\be
\label{sdform_1b}
\frac{\partial Q_i}{\partial t} = -\frac{1}{\Dx} 
\left(\apdqt_\imh +  \amdqt_\iph + \int_{x_\imh}^{x_\iph} A_i \, \qt_x dx\right).
\ee
%The fluctuations may be computed using a suitable (exact or approximate)
%Riemann solver.
In general, the integral in \eqref{sdform_1b} must be evaluated by quadrature; 
however, for the conservative system \eqref{conslaw}, the integral can be evaluated 
exactly, and is given by 
\be \label{fluxint}
\int_{x_\imh}^{x_\iph} A_i \, \qt_x dx = f(q^L_\iph)-f(q^R_\imh).
\ee
If the fluctuations are computed using a Roe solver or some other conservative
wave-propagation Riemann solver, then the flux difference appearing in
\eqref{fluxint} is equal to the sum of fluctuations from a fictitious ``internal''
Riemann problem for the current cell $i$, just as in the linear case above:
\be \label{fluxint-1}
f(q^L_\iph)-f(q^R_\imh) = \apdqt_i + \amdqt_i = \adqt_i.
\ee
Specifically, the fluctuations $\apmdqt_i$ are those resulting from the Riemann
problem with initial states $(q^R_\imh,q^L_\iph)$.
Then we can write \eqref{sdform_1b} also as
\be
\label{sdform_1c}
\frac{\partial Q_i}{\partial t} = -\frac{1}{\Dx} 
\left(\amdqt_\iph + \apdqt_\imh + \adqt_i\right).
\ee

Note that, for the conservative system \eqref{conslaw}, if a Roe solver or an
$f$-wave solver (see Section \ref{subsec:F-wave})
is used, then the fluctuations are equal to the flux differences
\begin{eqnarray}
\label{lfluct}
\amdqt_\imh & = & \hat{f}_\imh - f(q^L_\imh) \\
\apdqt_\imh & = & f(q^R_\imh) - \hat{f}_\imh,
\label{rfluct}
\end{eqnarray}
where $\hat{f}_\imh$ is the numerical flux at $x_\imh$.  Thus
\eqref{sdform_1c} is in that case equivalent
to the traditional flux-differencing method
\be
\dd{Q_i}{t} = -\frac{1}{\Dx}\left(\hat{f}_\iph-\hat{f}_\imh\right).
\ee
In particular, the scheme is conservative in this case.  

\subsection{Capacity-form differencing}
In many applications the system of conservation laws takes the form 
\be
\label{eq:capacitysys1D}
\kappa(x)q_t + f(q)_x = 0,
\ee
in one space dimension, or 
\be
\label{eq:capacitysys2D}
\kappa(x,y)q_t + f(q)_x + g(q)_y = 0,
\ee
in two dimensions, where $\kappa$ is a given function of space and is usually
indicated as \emph{capacity function} (see \cite{leveque1997}).  Systems like
(\ref{eq:capacitysys1D}) and (\ref{eq:capacitysys2D}) arise naturally in the
derivation of a conservation law, where the flux of a quantity is naturally
defined in terms of one variable $q$, whereas it is a different quantity
$\kappa q$ that is conserved. For instance, for the flow in a porous media,
$\kappa$ would be the porosity. Note that a capacity function can also appear
in systems that are not in conservation form, e.g. the quasilinear system
\begin{changebar}
(\ref{varlin}).
\end{changebar}

Several approaches can be used to reduce such a system to a more familiar conservation law. 
One natural approach is the capacity-form differencing \cite{leveque1997},
\be
\label{capacity-form-diff}
\frac{\partial Q_i}{\partial t} = -\frac{1}{\kappa_i \Dx} 
\left(\apdqt_\imh + \amdqt_\iph + \adqt_i\right),
\ee
where $\kappa_i$ is the capacity of the $i$th cell. This is a simple extension
of (\ref{sdform_1a}) or (\ref{sdform_1c}) which ensures that $\sum \kappa_i
Q_i$ is conserved (except possibly at the boundaries) and yet allows the
Riemann solution to be computed based on $q$ as in the case $\kappa = 1$.

\subsection{$f$-wave Riemann solvers and well-balancing\label{subsec:F-wave}}
For application to conservation laws, it is desirable to ensure that the
wave-propagation discretization is conservative.  This can easily be
accomplished by using
%In this section we review a technique for easily developing conservative
%wave propagation Riemann solvers for general nonlinear systems.
%These are 
an $f$-wave Riemann solver \cite{bale2002}.  Use of $f$-wave solvers is
also useful for problems with spatially varying flux function, 
as well as problems involving balance laws near steady state.
Further uses of the $f$-wave approach can be found in
\cite{ahmad:aiaa08,
ahmad-lindeman:fwave, khouider-majda:balanced, NormanNairSemazzi2010, sniv-pasko:fwave}.

The idea of the $f$-wave splitting for (\ref{conslaw}) is to decompose the flux
difference $f(q_r) - f(q_l)$ into waves rather than the $q$-difference used in
(\ref{eq:waveDecompFV}), i.e.  we decompose the flux difference as a linear
combination of the right eigenvectors $r^p$ of some Jacobian:
\be
\label{eqn:f-waveDef}
f(q_r) - f(q_l) = \sum_p \beta^p r^p = \sum_p \mathcal{Z}^p(q_l,q_r).
\ee
%The coefficients $\beta^p_\imh $ are given by
%\be
%\beta^p_\imh = l^p_\imh (f(Q_i) - f(Q_{i-1})),
%\ee
%where $l^p_\imh$ are the left eigenvectors of the Jacobian or approximate Jacobian matrix $A_\imh$. 
%Note that the $f$-waves have the dimensions of a $q$ increment multiplied by the wave speed.

The fluctuations are then defined as
\begin{align*}
\amdq & \equiv \sum_{p:s^p <0} \mathcal{Z}^p(q_l,q_r) &
\apdq & \equiv \sum_{p:s^p >0} \mathcal{Z}^p(q_l,q_r) &
\end{align*}
Note that the total fluctuation in cell $i$ is given simply by
$$
\adqt_i = f\left(q^L_\iph\right) - f\left(q^R_\imh\right).
$$

An advantage of particular interest is the possibility to include source terms directly into the $f$-wave decomposition. In fact, for balance laws that include non-hyperbolic terms,
$$
q_t + f(q)_x = \psi(q,x),
$$
one can easily extend this algorithm by first discretizing the source term to obtain values $\Psi_\imh$ at the cell interfaces and then compute the waves $\mathcal{Z}^p_\imh$ by splitting
\be \label{f-wave-balance}
f(q_r) - f(q_l) - \Dx \Psi(q_l,q_r,x) = \sum_p \beta^p r^p = \sum_p \mathcal{Z}^p(q_l,q_r).
\ee
%where
%\be
%\beta^p_\imh = l^p_\imh (f(Q_i) - f(Q_{i-1}) - \Dx \Psi_\imh).
%\ee
Here $\Psi(q_l,q_r,x)$ is some suitable average of $\psi(q,x),$ between the neighboring states.
In Bale et al. \cite{bale2002}, it has been shown that for the second-order FV
wave-propagation scheme implemented in Clawpack, the $f$-wave approach is very
useful for handling source terms, especially in cases where the solution is
close to a steady state because it leads to a well-balanced scheme.
However, for our high order wave-propagation scheme, application of
the $f$-wave algorithm with component-wise or characteristic-wise reconstruction \cite{Qiu2002187}
(which take no account of the source term) does not lead to a method that 
is well-balanced, even though the source term is accounted for in the Riemann
solves which compute $\apmdqt_\imh$. This is because with the aforementined WENO approaches the reconstructed solution within each cell is not constant (i.e. $q^L_\iph \neq q^R_\imh$) and $\adqt_i \neq 0$.
%The total fluctuation is computed by solving a Riemann problem with left state $q^+_\imh$ and right state $q^-_\iph$ which is decoupled from the Riemann problems at the two cell's interfaces. Indeed, the fluctuations and the total fluctuation are computed using different Jacobian matrices, i.e. different eigenstructures. 
%We remark that in general it may be difficult to find a suitable average $\Psi$
%that leads to well-balancing \cite{rjl:wbfwave10}.  
%However, in the case of the steady-state-at-rest
%for the shallow water equations it is straightforward, as demonstrated in Section
%\ref{subsubsec:2Dperturb-SWEs}.

In this work, we consider an extension of the $f$-wave well-balancing
technique that is compatible with our higher order scheme.  The extension
is useful for problems in which the source term vanishes over the interior
of each cell (i.e., its effect is concentrated at cell interfaces).  This
is the case, for instance, when considering the shallow water equations with
piecewise-constant bathymetry (see Section \ref{subsubsec:2Dperturb-SWEs}).
This well-balancing technique uses the $f$-wave Riemann solver in combination 
with a $f$-wave-slope reconstruction (see Section 
\ref{subsect:f-wave-slope-recon}). With this approach, the contribution of the source 
term is directly included in the Riemann solve at the cell's interface, 
i.e. (\ref{f-wave-balance}).
The reconstruction methods considered in this work are presented in the next section.

%Assuming that the source term $\psi$ is just a function of the spatial 
%coordinates (e.g. geometric source term in the shallow water equations) and 
%it is continuos at the cell faces, one approach to achieve high-order 
%well-balancing is taking into account 
%the contribution of $\psi$ by using the $f$-wave algorithm in the calculation 
%of the total fluctuations (\ref{eq:3intsc-fluct}). 
%This requires an interpolation of the source term at the Gauss 
%quadrature points for an high-order accurate evaluation of its 
%integral over the cell volume. The interpolation procedure is the key step to 
%achieve high-order well-balancing. It is very challenging to find a
%multidimensional algorithm that interpolates with high-order accuracy a general
%source term and leads to a non-oscillatory solution which cancels out the 
%total fluctuation contribution. Therefore, in this work we propose an 
%alternative approach which is just second-order 
%accurate in space but it can be easily applied to any balance law with spatially
%varying source term.
%%%%
\subsection{Reconstruction\label{subsec:waveSlopeReconstruction}}
The reconstruction \eqref{eq:pwp} should be performed in a manner that yields high order
accuracy but avoids spurious oscillations near discontinuities.  For this purpose,
we use weighted essentially non-oscillatory (WENO) reconstruction \cite{shu2009}.
The spatial accuracy of the method will in general be equal to that of the reconstruction.
In the present work we employ fifth-order WENO reconstruction.

\subsubsection{Scalar WENO Reconstruction}
WENO reconstruction formulas are typically written in terms of the divided differences
$\DQ_{\ipmh}/\Dx$.  It is possible to rewrite them in terms of the difference ratios
\be \label{theta-scalar}
\theta_{\imh,j} = \frac{\DQ_{\imh+j}}{\DQ_\imh}.
\ee
as long as $\DQ_\imh\ne0$.  Then the reconstructed interface values in
cell $i$ are given by
\begin{subequations} \label{recon-scalar}
\begin{align} 
q^R_\imh & = Q_i - \phi(\theta_{\imh,2-k},\dots,\theta_{\imh,k-1})\DQ_\imh \\
         %= Q_i - \phi_\imh \DQ_\imh \\
q^L_\iph & = Q_i + \phi(\theta_{\iph,1-k},\dots,\theta_{\iph,k-2})\DQ_\imh,
         %= Q_i + \phi_\iph \DQ_\iph
\end{align}
\end{subequations}
where $\phi$ is a particular nonlinear function that we will not write out here.
The usefulness of \eqref{recon-scalar} is that it allows WENO reconstruction to
be applied to waves directly by redefining $\theta$, as we do below.
In the case that $\DQ_\imh\approx0$ (to near machine precision), the value of
$\phi$ is set to zero.

For systems of equations, the simplest approach to reconstruction is
component-wise reconstruction, which consists of applying the
scalar reconstruction approach \eqref{theta-scalar}-\eqref{recon-scalar} to
each element of $q$.  A more sophisticated
approach is characteristic-wise reconstruction, in which an eigendecomposition
of $q$ is performed, followed by reconstruction of each eigencomponent.
For problems with spatially-varying coefficients, even the characteristic-wise
reconstruction may not be satisfying, since it involves comparing coefficients of
eigenvectors whose direction in state space varies from one cell to the next.
In Clawpack, an alternative kind of TVD limiting known as {\em wave limiting} 
has been implemented and shown to be effective for such problems.

\subsubsection{Wave-slope reconstruction}
In order to implement a wave-type WENO limiter, we first solve a Riemann problem
at each interface $x_\imh$, using the adjacent cell average 
values $Q_{i-1},Q_i$ as left and right states.
This results in a set of waves $\Wop^p_\imh$, which are used solely for the purpose
of reconstruction.  This reconstruction is performed by replacing \eqref{theta-scalar}
by
\begin{align} \label{theta-wave}
\theta^p_{\imh,j} & = \frac{\Wop^p_{i-\frac{1}{2}+j}\cdot\Wop^p_\imh}{\Wop^p_\imh\cdot\Wop^p_\imh}
\end{align}
and replacing \eqref{recon-scalar} by
\begin{subequations} \label{recon-wave}
\begin{align}
q^L_\imh & = Q_{i-1} + \sum_p \phi(\theta^p_{\iph,1-k},\dots,\theta^p_{\iph,k-2}) \Wop^p_\imh \\
q^R_\imh & = Q_i \ \ \ - \sum_p \phi(\theta^p_{\imh,2-k},\dots,\theta^p_{\imh,k-1}) \Wop^p_\imh,
\end{align}
\end{subequations}
This approach takes into account the smoothness of the $p$th 
characteristic component of the solution by using the information arising from 
the $k$-cell stencils. It is intended to be similar to that used in Clawpack 
\cite{levequefvbook} and can be conveniently implemented using the same Riemann
solvers supplied with Clawpack.

\subsubsection{$f$-wave-slope reconstruction}\label{subsect:f-wave-slope-recon}
Wave-slope reconstruction can also be performed using an $f$-wave Riemann solver.
This is useful for computing near-equilibrium solutions of balance laws.  The procedure
is identical to that above, except that the Riemann problem is solved with the $f$-wave Riemann solver at each interface
$x_\imh$, using the adjacent cell average values $Q_{i-1},Q_i$ as left and
right states.
Since the resulting $f$-waves have the form of a $q$ increment multiplied by the wave
speed \cite{bale2002}, the waves $\Zop$ are normalized by the corresponding wave speed before being
used for reconstruction:
\begin{align} \label{wavediv}
\Wop^p_\imh & = \mathcal{Z}^p_\imh / s^p_\imh.
\end{align}
%Particular attention must be given to the special situations of $s^p = 0$ which 
%requires $\Wop^p=0$ when the wave-propagation approach leads to a conservative
%algorithm \cite{bale2002}.
In this paper we assumed that the original hyperbolic system has no zero 
eigenvalues ($s^p \neq 0$) or eigenvalues that change sign between grid cells 
(i.e. the resonant case). The reconstruction procedure \eqref{theta-wave}-\eqref{recon-wave} is then applied
to the waves computed by \eqref{wavediv}.

For systems with source terms that are at a steady state, assuming that
\eqref{f-wave-balance} holds, results in $\Zop^p=0$
for all waves. Therefore, the $f$-wave-slope reconstruction
will yield to a constant reconstructed function in regions
where source terms and hyperbolic terms are balanced.  Then all fluctuations
computed in the update step will vanish, so the steady state will be preserved
exactly.

\subsection{Time integration}
The semi-discrete scheme can be integrated in time using any
initial-value ODE solver. Herein we use the ten-stage fourth-order 
strong-stability-preserving Runge-Kutta scheme of \cite{ketcheson2008}. 
This method has a large stability region and a large SSP coefficient, allowing
use of a relatively large CFL number in practical computations. In all numerical examples
of the next section, a CFL number of 2.45 is used.

To summarize, the full semi-discrete algorithm used in each Runge-Kutta stage
is as follows.
\begin{enumerate}
\setcounter{enumi}{-1}
  \item {\em (only if using wave-slope reconstruction)} Solve the Riemann problem at
        each interface $x_\imh$  using the adjacent cell average  values $Q_{i-1},Q_i$ 
        as left and right states. %This step results in a set of waves $\Wop^p_\imh$.
  \item Compute the reconstructed piecewise function $\qt$, and in particular the
        states $q^R_\imh,q^L_\iph$ in each cell, using 
        component-wise, characteristic-wise, or wave-slope reconstruction.
  \item At each interface $x_\imh$, compute the fluctuations $\apdqt_\imh$ and $\amdqt_\imh$
        by solving the Riemann problem with initial states $(q^L_\imh,q^R_\imh)$.
  \item Over each cell, compute the integral $\int A\,\qt_x dx$.  For conservative
        systems this is just the total fluctuation $\adq_i$.
  \item Compute $\partial Q/\partial t$ using the semi-discrete scheme \eqref{sdform_1b}.
%\item Using the wave speeds $s^p_\imh$ the CFL condition is inspected.
%\begin{itemize}
%\item If it is satisfied the solution at the new time level is computed by adding the cell averages updates obtained from the ODE solver.
%\item If it is not satisfied a smaller time step $\Delta t$, based on wave speeds $s^p_\imh$, is taken in the ODE solver and a new residual is computed.
%\end{itemize}
\end{enumerate}
Note again that, for conservative systems, the quadrature in step 3
requires nothing more than evaluating and differencing the fluxes.

%%%%

\subsection{Extension to Two Dimensions}
In this section, we extend the numerical wave propagation method
to two dimensions using a simple dimension-by-dimension approach.
The method is applicable to systems of the form
\be
q_t + A(x,y)q_x + B(x,y)q_y = 0
\ee
on uniform Cartesian grids.  

The 2D analog of the semi-discrete scheme \eqref{sdform_1c} is
\begin{align}
\label{sdform_2D}
\begin{split}
\frac{\partial Q_{ij}}{\partial t} = -\frac{1}{\Dx\Dy} 
 \left(\Aop^-\dqt_{\iph,j} + \Aop^+\dqt_{\imh,j} + \Aop\dqt_{i,j} \right. \\
  \left. +\Bop^-\dqt_{i,\jph} + \Bop^+\dqt_{i,\jmh} + \Bop\dqt_{i,j} \right).
\end{split}
\end{align}
For the method to be high order accurate, the fluctuation terms
like $\Aop^-\dqt_{\iph,j}$ should involve integrals over cell edges,
while the total fluctuation terms like $\Aop\dqt_{i,j}$ should involve
integrals over cell areas.  This can be achieved by forming a 
genuinely multidimensional reconstruction of $q$ and using, e.g.,
Gauss quadrature. An implementation following this approach 
exists in the SharpClaw software. For nonlinear problems containing shocks, 
the genuinely multidimensional reconstruction has been found
to be inefficient (at least for some simple test problems), as it typically 
yields only a small improvement in accuracy over the dimension-by-dimension 
scheme given below (which formally only second-order accurate), but has a much
greater computational cost on the same
mesh. Similar observations have been reported in Zhang et al. \cite{Zhang-2011-WENO},
where both approaches have been throughly tested and compared for linear and 
nonlinear systems. For problems with shocks, at least for the simple test 
problems presented, the two schemes give comparable resolution on 
the same meshes, despite their difference in formal order of accuracy. 
%it has been found that the dimension-by-dimension scheme is only 
%second-order accurate for nonlinear systems, whereas the genuinely 
%multidimensional reconstruction is still high order accurate. 

We now describe the dimension-by-dimension scheme for a single 
Runge-Kutta stage.  We first reconstruct piecewise-polynomial 
functions $\qt_j(x)$ along each row of the grid and $\qt_i(y)$ 
along each column, by applying a 1D
reconstruction procedure to each slice.  We thus obtain reconstructed
values
\begin{subequations}
\begin{align}
\qt_j^R(x_\imh) & \approx q(x_\imh,y_j) \\
\qt_j^L(x_\iph) & \approx q(x_\iph,y_j) \\
\qt_i^R(y_\imh) & \approx q(x_i,y_\imh) \\
\qt_i^L(y_\iph) & \approx q(x_i,y_\iph)
\end{align}
\end{subequations}
for each cell $i,j$.  The fluctuation terms in \eqref{sdform_2D}
are determined by solving Riemann problems between the appropriate
reconstructed values; for instance $\Bop^-\dqt_{i,\jph}$ is determined
by solving a Riemann problem in the $y$-direction with initial states
$(q^L_{i,\jph},q^R_{i,\jph}).$
In the case of conservative systems or piecewise-constant
coefficients, the total fluctuation terms $\Aop\dqt_{i,j}$ and
$\Bop\dqt_{i,j}$ can be similarly determined by summing the left- and right-going
fluctuations of an appropriate Riemann problem.  Thus, for instance,
%$\Bop\Dq_{i,j}$ is determined by the fluctuations resulting from the
%Riemann problem with initial data
%\begin{equation*}
%q(x,y,0)=\left\{ \begin{array}{cc} q^R_{i,\jmh} & x=x_i, \quad y<y_i \\ q^L_{i,\jph} & x=x_i, \quad y>y_i. \end{array}\right.
%\end{equation*}
$\Bop\dqt_{i,j}$ is determined by solving
Riemann problem in the $y$-direction with initial states $(q^R_{i,\jmh},q^L_{i,\jph}).$

%==================================================
\section{Numerical applications} \label{sect:applications}
%==================================================

In this section we present results of numerical tests using the wave
propagation methods just described.  The examples included are chosen to
emphasize the advantages of the wave propagation approach.
We make some comparisons with the well-known second-order wave propagation code
Clawpack \cite{clawpack45,levequefvbook}. 

%==================================================
\subsection{Acoustics} \label{sect:acoustics}
%==================================================
In this section, the high-order wave propagation algorithm is applied
to the 1D equations of linear acoustics in piecewise homogeneous materials:
\begin{subequations} \label{eq:acoustics}
\begin{eqnarray}
    p_t + K(x,y) (u_x + v_y) & = & 0 \\
    u_t + \frac{1}{\rho(x,y)} p_x & = & 0 \\
    v_t + \frac{1}{\rho(x,y)} p_y & = & 0.
\end{eqnarray}
\end{subequations}
Here $p$ is the pressure and $u,v$ are the x- and y-velocities, respectively.
The coefficients $\rho$ and $K$, which vary in space, are the density and bulk
modulus of the medium.  We will also refer to the sound speed $c=\sqrt{K/\rho}$.
Notice that in general since $\rho$ varies in space, the last two equations above
are not in conservation form.
This test case demonstrates that the proposed approach is able to 
solve hyperbolic system of equations written in nonconservative form. 
Of course, this system can be written in conservative form as follows:
\begin{subequations} \label{eq:acoustics_cons}
\begin{eqnarray}
    \epsilon_t - (u_x + v_y) & = & 0 \\
    \rho(x,y) u_t - (K(x,y)\epsilon)_x & = & 0 \\
    \rho(x,y) v_t - (K(x,y)\epsilon)_y & = & 0,
\end{eqnarray}
\end{subequations}
Where $\epsilon=-p/K$ is the strain.
As we will see, the latter form may be advantageous in terms of the accuracy
that can be obtained.

The grid is chosen so that the material is homogeneous in each computational cell,
and an exact Riemann solver is used at each interface; for details of this solver see e.g.
\cite{fogarty1999}.
%, so that
%in cell $i$ the Jacobian matrix is
%\begin{align}
%    %q & = \left( \begin{array}{c} p \\ u \end{array} \right), &
%    \mA_i & = \left( \begin{array}{cc}
%        0 & K_i \\
%        1/\rho_i & 0
%    \end{array} \right).
%\end{align}
%For a detailed description of the eigen-structure of this matrix we refer to \cite{leveque1997}.

%Let $Z_i=\sqrt{K_i\rho_i}$ denote the impedance and $c_i=\sqrt{K_i/\rho_i}$ the
%sound speed.  Then the eigenvectors of $\mA_i$ are
%\begin{align}
%r^1_i & = \left(\begin{array}{c} -Z_i \\ 1 \end{array}\right), &
%r^2_i & = \left(\begin{array}{c}  Z_i \\ 1 \end{array}\right),
%\end{align}
%and the eigenvalues are
%\be
%\lambda^1_i=-c_i, \quad \quad \lambda^2_i=c_i.
%\ee
%At an interface between different materials, the Riemann solution consists
%of two waves that are multiples of the appropriate eigenvector in each material:
%\begin{align}
%\Wop^1_\imh & = \alpha^1_\imh\left(\begin{array}{c} -Z_{i-1} \\ 1 \end{array}\right), &
%\Wop^2_\imh & = \alpha^2_\imh\left(\begin{array}{c}  Z_i     \\ 1 \end{array}\right).
%\end{align}
%These waves move at speeds
%\be
%\lambda^1_\imh=-c_{i-1}, \quad \quad \lambda^2_\imh=c_i,
%\ee
%respectively.  The wave strengths are 
%\begin{subequations}
%\begin{align}
%\alpha^1_\imh & = \frac{-\Delta p_\imh + Z_i    \Delta u_\imh}{Z_i+Z_{i-1}} \\
%\alpha^2_\imh & = \frac{ \Delta p_\imh + Z_{i-1}\Delta u_\imh}{Z_i+Z_{i-1}}.
%\end{align}
%\end{subequations}

\subsubsection{One-dimensional acoustics}
We first consider one-dimensional acoustic waves in a piecewise-constant medium
with a single interface.  Namely, we solve \eqref{eq:acoustics} on the interval 
$x\in[-10,10]$ with
$$
(\rho,c) = \left\{\begin{array}{cr}
    (\rho_l,c_l) & x<0 \\
    (\rho_r,c_r) & x>0 \end{array}\right.
$$
We measure the convergence rate of the solution
in order to verify the order of accuracy for smooth solutions.  
%We thus require an initial condition that is $p$-times differentiable,
%where $p$ is greater than the expected order of accuracy.  
%Additionally, we require that the initial condition have compact support
%so that after some time the solution will be around the interface.
The initial condition is a compact, six-times differentiable purely right-moving 
pulse:  %given by a Newton interpolating polynomial,
\begin{align*} %\label{newtonIC}
p(x,0) & = \frac{((x-x_0)-a)^6((x-x_0)+a)^6}{a^{12}} \xi(x-x_0) \\
u(x,0) & = p(x,0)/Z(x)
\end{align*}
where
$$
\xi(x-x_0)=\left\{ \begin{array}{cc} 
0 & \mbox{for } |x-x_0|> a) \\
1 & \mbox{for } |x-x_0| \le a.
\end{array} \right.
$$
with $x_0=-4$ and $a=1$, and $Z(x) = \sqrt{K(x)\rho(x)}$.
This condition was chosen to be sufficiently smooth to demonstrate the
design order of the scheme and to give a solution that is identically zero at the
material interface at the initial and final times.

Table \ref{tbl:homogeneous} shows $L_1$ errors and convergence
rates for propagation in a homogeneous medium with $\rho_l=c_l=\rho_r=c_r=1$. 
Here we use componentwise reconstruction.  Specifically, we compute
\be
\label{eq:defL1Err}
E_{L_1} = \Dx\sum_i |Q_i - \bar{Q}_i|
\ee
where $\bar{Q}_i$ is a highly accurate solution cell average computed by characteristics
or by using a very fine grid.  For the acoustics problems in this section,
$\bar{Q}$ is computed using characteristics and adaptive Gauss quadrature.
Table \ref{tbl:homogeneous} indicates that in each case, the order
of convergence is approximately equal to the design order of the discretization.

%  Errors computed by clawtest.py on Tue May 26 2009 at 22:50:37 

      \begin{table}[th]
      \caption{Errors for homogeneous acoustics test \label{tbl:homogeneous}}
      \begin{center}
    \begin{tabular}{c|cc|cc}
& \multicolumn{2}{c|}{SharpClaw}& \multicolumn{2}{c}{Clawpack}\\ \hline 
mx & Error & Order& Error & Order\\ \hline 
 200 &         3.60e-02 &              &         4.10e-02 &              \\ 
\shadeRow 400 &         3.65e-03 &             3.30&         1.30e-02 &             1.66\\ 
 800 &         1.85e-04 &             4.31&         3.61e-03 &             1.85\\ 
\shadeRow1600 &         7.35e-06 &             4.65&         8.94e-04 &             2.01\\ 
%3200 &         2.72e-07 &             4.76&         2.19e-04 &             2.03\\ 

      \hline
      \end{tabular}
      \end{center}
      \end{table}

To test the accuracy in the presence of discontinuous coefficients we take
$$\rho_l=c_l=1 \quad \rho_r=4 \quad c_r=1/2,$$
with an impedance ratio of $Z_r/Z_l=2$.
As was noted in \cite{bale2002}, this system can also be solved
in the conservative form \eqref{eq:acoustics_cons} using the $f$-wave approach.  We include results
of this approach, where we have also performed characteristic-wise
rather than component-wise reconstruction.
Results are shown in Table \ref{tbl:interface 1}.
In this case all schemes exhibit a convergence
rate below the formal order, even though the initial
and final solutions are smooth.  To investigate this further,
we repeat the same test with a wider pulse by taking $a=4$.
Results are shown in Table \ref{tbl:interface 1 wide}.

For the latter test, we observe a convergence rate of approximately
two for SharpClaw, one for Clawpack, and five for SharpClaw using the $f$-wave
approach and characteristic-wise reconstruction.  The last convergence
rate is remarkable, considering that the solution is not differentiable
when it passes through the material interface.  Further investigation
of the accuracy of this approach for more complicated problems
with discontinuous coefficients is ongoing.  
In tests not shown here, Clawpack
achieves approximately second-order accuracy when used with an
$f$-wave Riemann solver for this problem.
%For this problem, the errors can be divided into two sources.
%First, the truncation error that occurs as the pulse is
%propagated in a homogeneous medium; this error is of the
%order given by the design order of each method.  This error
%is large for a narrow pulse that is not well resolved on the
%grid.

%The second source of error arises in the reconstruction step, 
%when the solution is reconstructed using stencils that cross the
%material interface.  The high-order reconstruction is based 
%on an assumption of smoothness of the solution, which does
%not hold at the interface.  Since the jump in the first
%derivative of the solution at the interface is $\Oop(1)$,
%the error in the reconstructed values in cells whose stencil
%overlaps the interface is $\Oop(\Dx)$.  Since the total area
%of all such cells is $\Oop(\Dx)$, the resulting global 
%error is $\Oop(\Dx^2)$.
%This error dominates when the pulse is well-resolved on the
%grid, so that the first type of error is small.

\begin{table}[th]
      \caption{Errors for acoustics interface with narrow pulse \label{tbl:interface 1}}
      \begin{center}
    \begin{tabular}{c|cc|cc|cc}
& \multicolumn{2}{c|}{SharpClaw}& \multicolumn{2}{c|}{SharpClaw $f$-wave}& \multicolumn{2}{c}{Clawpack}\\ \hline 
mx & Error & Order& Error & Order& Error & Order\\ \hline 
 200 &         2.10e-01 & & 9.50e-02 &              &         1.98e-01 &              \\ 
\shadeRow 400 &         5.98e-02 &             1.81& 1.42e-02 & 2.74 &        7.26e-02 &             1.45\\ 
 800 &         1.25e-02 &             2.26& 1.42 e-03 & 3.32 &        2.21e-02 &             1.71\\ 
\shadeRow1600 &         1.17e-03 &             3.42& 1.20e-04 & 3.56 &         7.86e-03 &             1.49\\ 
      \hline
      \end{tabular}
      \end{center}
      \end{table}

%  Errors computed by clawtest.py on Wed May 27 2009 at 10:31:51 

      \begin{table}[th]
      \caption{$L^1$ Errors for acoustics interface problem with wide pulse (a=4)\label{tbl:interface 1 wide}}
      \begin{center}
    \begin{tabular}{c|cc|cc|cc}
& \multicolumn{2}{c|}{SharpClaw}& \multicolumn{2}{c|}{SharpClaw $f$-wave}& \multicolumn{2}{c}{Clawpack}\\ \hline 
mx & Error & Order& Error & Order & Error & Order\\ \hline 
 200 &         9.67e-03 & &  5.01e-03 &             &         5.23e-02 &              \\ 
\shadeRow 400 &         2.01e-03 &             2.27& 4.63e-04 & 3.44 &         2.32e-02 &             1.17\\ 
 800 &         4.89e-04 &             2.04& 2.51e-05 & 4.36 &        1.09e-02 &             1.09\\ 
\shadeRow1600 &         1.22e-04 &             2.00& 6.49e-07 & 5.12 &         5.26e-03 &             1.05\\ 
      \hline
      \end{tabular}
      \end{center}
      \end{table}

\subsubsection{A Two-dimensional sonic crystal}
In this section we model sound propagation in a sonic crystal.
A sonic crystal is a periodic structure composed of materials
with different sounds speeds and impedances.  The periodic
inhomogeneity can give rise to {\em bandgaps} -- frequency bands
that are completely reflected by the crystal.  This phenomenon
is widely utilized in photonics, but its significance for acoustics
has only recently been considered.  Photonic crystals can be
analyzed quite accurately using analytic techniques, since they
are essentially infinite size structures relative to the wavelength
of the waves of interest.  In contrast,
sonic crystals are typically only a few wavelengths in size,
so that the effects of their finite size cannot be neglected.
For more information on sonic crystals, see for instance the
review paper \cite{miyashita2005}.

We consider a square array of square rods in air with a plane wave 
disturbance incident parallel to one of the axes of symmetry.  The array 
is infinitely wide but only eight periods deep.  The lattice spacing is 10 cm
and the rods have a cross-sectional side length of 4 cm, so that the filling
fraction is $0.16$.  This crystal is similar to one studied in
\cite{sanchis2001}, and it is expected that sound waves in the 1200-1800 Hz
range will experience severe attenuation in passing through it, while
longer wavelengths will not be significantly attenuated.

A numerical instability very similar to that observed in 1D
simulations in \cite{fogarty1998,fogarty1999} was observed 
when the standard Clawpack method was applied to this problem.  
The fifth-order  WENO method with characteristic-wise limiting showed no 
such instability.

\begin{figure}
\centering
\includegraphics[width=4in]{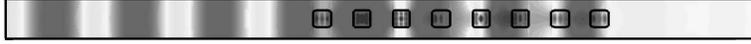}
\caption{Pressure in the sonic crystal for a long wavelength plane wave incident from the left.\label{sclongwave}}
\end{figure}

\begin{figure}
\centering
\includegraphics[width=3in]{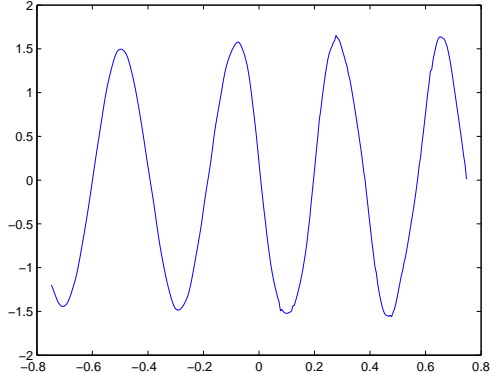}
\caption{Pressure in the sonic crystal for a long wavelength plane wave incident from the left.\label{sclongwave_slice}}
\end{figure}

Figure \ref{sclongwave} shows a snapshot of the RMS pressure distribution
in space for a plane wave with $k=15$ incident from the left.  
The RMS (root mean square) pressure is computed as follows:
%at each point by taking
%the square of the pressure, averaging over one temporal period of the plane
%wave, and taking the square root:
\begin{align}
p_\textup{RMS}(x,y) = \sqrt{\frac{1}{T}\int_t^{t+T}p^2(x,y,t)dt}.
\end{align}
This wave has a
frequency of about 800 Hz, well below the partial
band gap.  As expected, the wave passes through the crystal without
significant attenuation.  In Figure \ref{sclongwave_slice}, the pressure
is plotted along a slice in the $x$-direction approximately midway 
between rows of rods.

\begin{figure}
\centering
\includegraphics[width=4in]{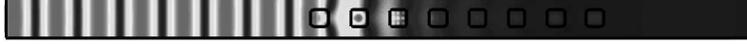}
\caption{Snapshot of RMS pressure distribution in space in the sonic crystal for a plane wave incident from the left.\label{sc1}}
\end{figure}

\begin{figure}
\centering
\includegraphics[width=3in]{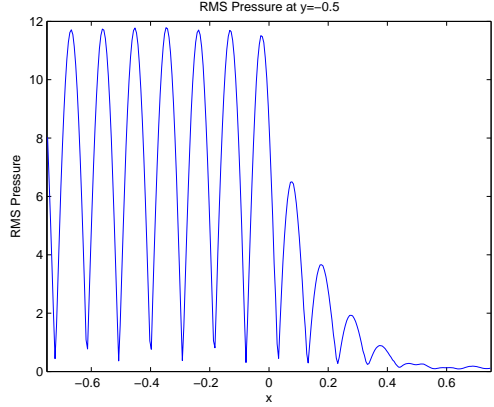}
\caption{Snapshot of RMS pressure distribution in space in the sonic crystal along a slice.\label{sc1_slice}}
\end{figure}

Figure \ref{sc1} shows a snapshot of the RMS pressure distribution in space
for an incident plane wave with with frequency 1600 Hz, inside the partial bandgap.
%$k = 29.22 \mbox{ m}^{-1}, c=344$ m/s.
Notice that the wave is almost entirely reflected, resulting in a 
standing wave in front of the crystal.  Figure \ref{sc1_slice} 
shows the RMS pressure along a slice in the $x$-direction.

\subsection{Nonlinear Elasticity in a Spatially Varying Medium}
\label{subsec:stego}
In this section we consider a more difficult test involving nonlinear
wave propagation and many material interfaces.  This problem was considered
previously in \cite{leveque2002} and studied extensively in \cite{leveque2003}.
Solitary waves were observed to arise from the interaction
of nonlinearity and an effective dispersion due to material
interfaces in layered media.

Elastic compression waves in one dimension are governed by the equations
\begin{subequations} \label{nel_pde}
\begin{align}
\epsilon_t(x,t)-u_x(x,t) & = 0 \\
(\rho(x)u(x,t))_t - \sigma(\epsilon(x,t),x)_x & = 0.
\end{align}
\end{subequations}
where $\epsilon$ is the strain, $u$ the velocity,
$\rho$ the density, and $\sigma$ the stress.
This is a conservation law of the form \eqref{conslaw}, with
\begin{align*}
q(x,t) & = \left(\begin{array}{c} \epsilon \\ \rho(x) u \end{array}\right)
& f(q,x) & = \left(\begin{array}{c} -u \\ -\sigma(\epsilon,x) \end{array}\right).
\end{align*}
Note that the density and the stress-strain relationship vary
in $x$.  %We will also refer to the sound speed $c(x)$, impedance $Z(x)$,
%and linearized bulk modulus $K(x)$, given by
%\begin{align}
%c(x) & = \sqrt{\sigma_\epsilon(\epsilon,x)/\rho(x)} \\
%Z(x) & = \rho(x) c(x) \\
%K(x) & = \left.\sigma_\epsilon(\epsilon,x)\right|_{\epsilon=0}.
%\end{align}
The Jacobian of the flux function is
\begin{align*}
f'(q) = \left(\begin{array}{cc} 0 & -1/\rho(x) \\ -\sigma_\epsilon(\epsilon,x) & 0 \end{array}\right).
\end{align*}
%with eigenvectors
%\begin{align}
%r^1 & = \left(\begin{array}{c} 1 \\-Z(\epsilon,x) \end{array}\right), &
%r^2 & = \left(\begin{array}{c} 1 \\ Z(\epsilon,x) \end{array}\right)
%\end{align}

In the case of the linear stress-strain relation 
$\sigma(x)=K(x)\epsilon(x)$, \eqref{nel_pde} 
is equivalent to the one-dimensional form of the acoustics equations 
considered in the previous section.

%\subsubsection{An F-wave Riemann Solver}
%To apply wave propagation methods to 
%\eqref{nel_pde}, we need to define a Riemann solver for this system.
%We will use the $f$-wave solver developed in \cite{leveque2002}.
%We assume that $\rho$ is constant in each cell and that $\sigma$ depends
%only on $\epsilon$ (not explicitly on $x$) within a given cell; we denote their
%values in cell $i$ by $\rho_i,\sigma_i(\epsilon)$.
%Given a Riemann problem at $x_\imh$ we use the approximate wave speeds
%\begin{align}
%s^1_\imh & = -\sqrt{\frac{\sigma'_{i-1}(\epsilon^-_\imh)}{\rho^-_\imh}} &
%s^2_\imh & =  \sqrt{\frac{\sigma'_{i}  (\epsilon^+_\imh)}{\rho^+_\imh}}.
%\end{align}
%and the eigenvectors
%\begin{align}
%r^1_\imh & = \left(\begin{array}{c} 1 \\-\sqrt{\rho^-_\imh \sigma'_{i-1}(\epsilon^-_\imh)} \end{array}\right), &
%r^2_\imh & = \left(\begin{array}{c} 1 \\ \sqrt{\rho^+_\imh \sigma'_{i}(\epsilon^+_\imh)} \end{array}\right).
%\end{align}
%
%Using (\ref{eqn:f-waveDef}), the flux difference is then decomposed as
%\be
%f(q^+_\imh)-f(q^-_\imh) = \beta^1_\imh r^1_\imh + \beta^2_\imh r^2_\imh.
%\ee
%Then the fluctuations are simply
%\begin{align}
%\Aop^-\Dq_\imh & = \beta^1_\imh r^1_\imh &
%\Aop^+\Dq_\imh & = \beta^2_\imh r^2_\imh &
%\end{align}

We consider the piecewise constant medium studied in 
\cite{leveque2002,leveque2003}:
\begin{align} \label{eq:layered}
(\rho(x),K(x)) & = \left\{\begin{array}{cl}
    (1,1) & \mbox{if } j<x<(j+\frac{1}{2})
        \mbox{ for some integer } j \\
    (4,4) & \mbox{otherwise},
    \end{array}\right.
\end{align}
with exponential stress-strain relation
\be \label{eq:exp-stress} \sigma(\epsilon,x)=\exp(K(x)\epsilon)-1. \ee
The initial condition is uniformly zero, and the boundary
condition at the left generates a half-cosine pulse.

\begin{figure}
\centerline{
\subfigure[Strain.\label{fig:stego_comp_strain}]{\includegraphics[width=3in]{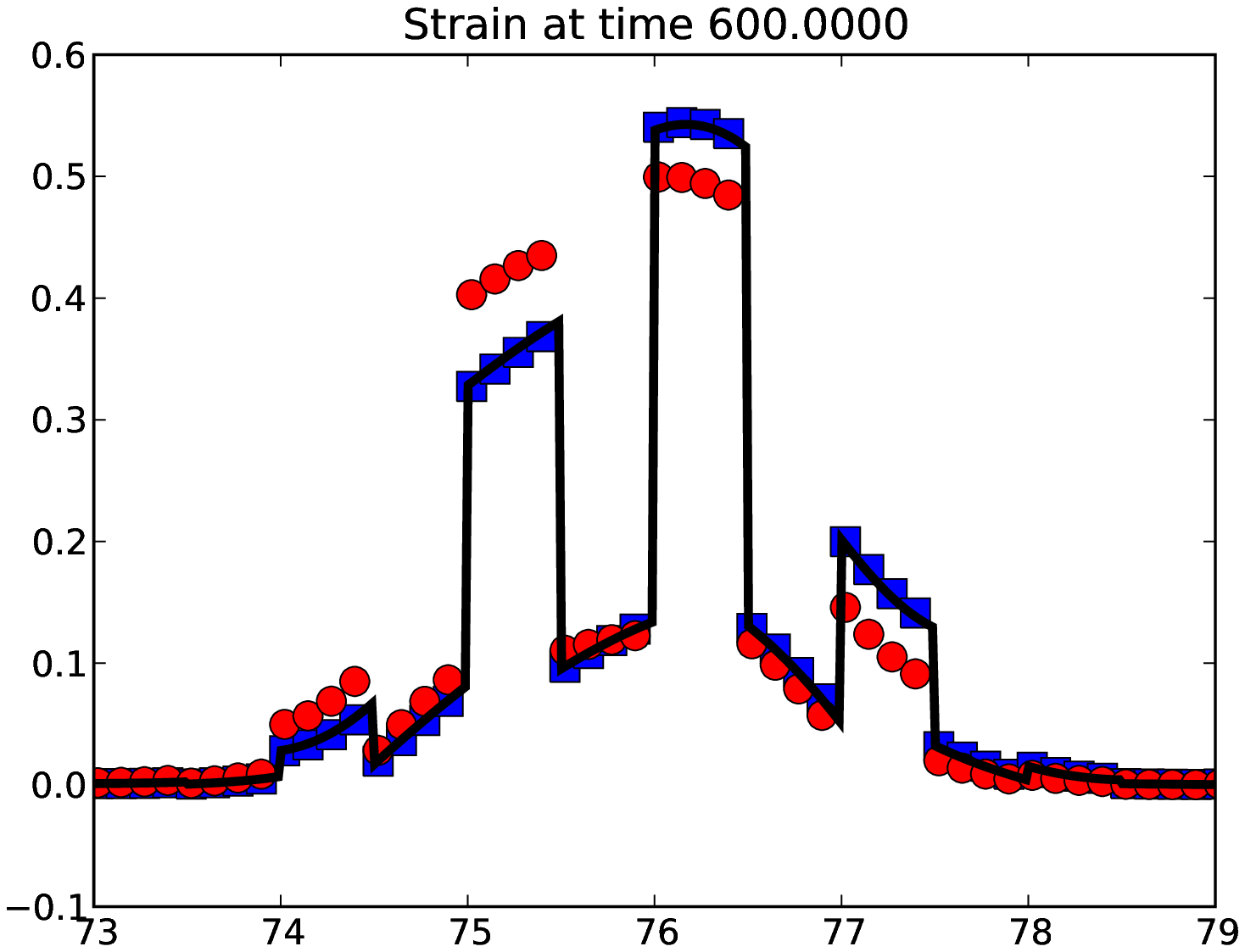}}
\subfigure[Stress.\label{fig:stego_comp_stress}]{\includegraphics[width=3in]{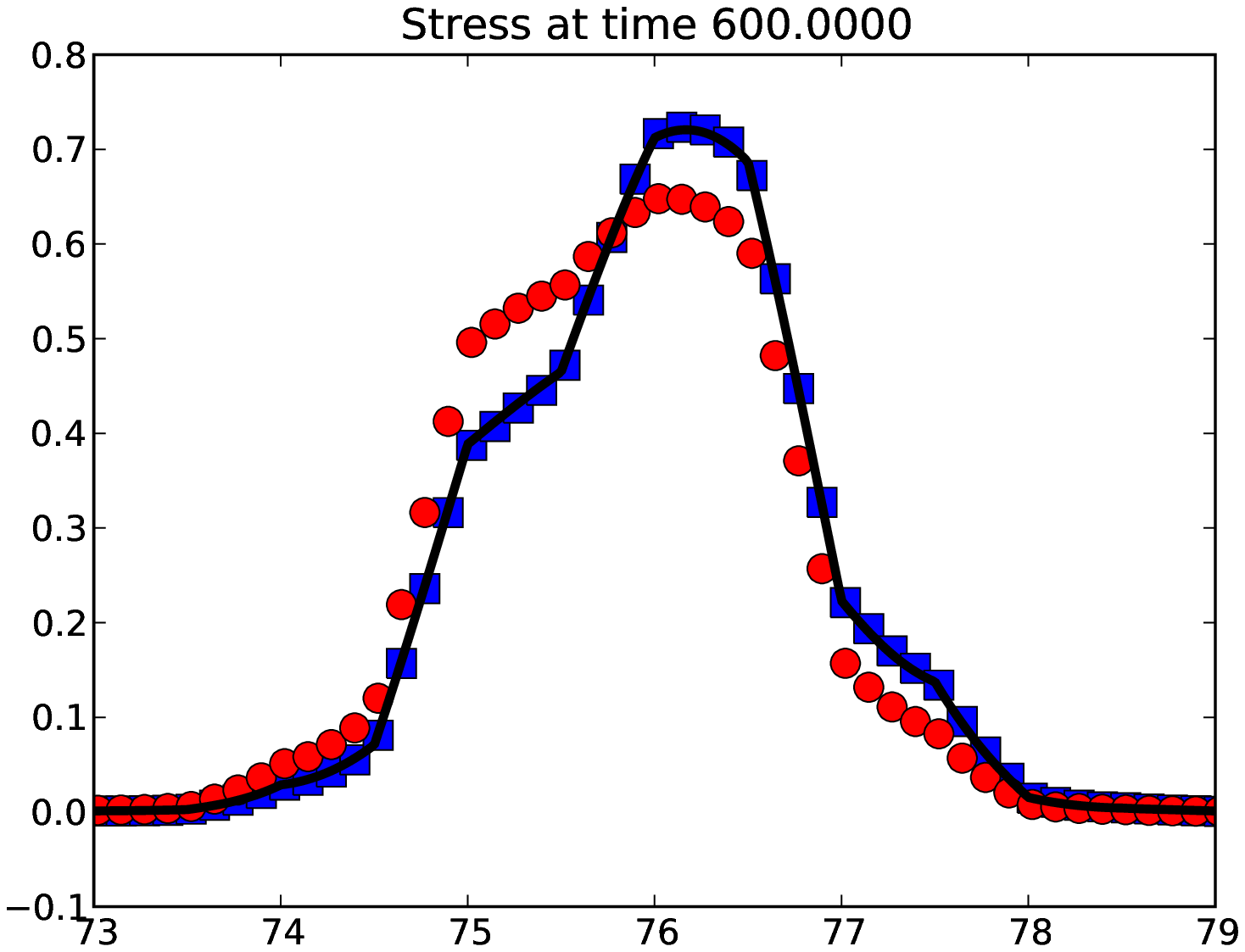}}}
\caption[Comparison of Clawpack and SharpClaw solutions of Stegoton problem.]
{Comparison of Clawpack (red circles) and SharpClaw (blue squares)
solution of the stegoton problem using 24 cells per layer.  For clarity, 
only every third solution point is plotted.  The black line represents a 
very highly resolved solution.\label{fig:stego_comp}}
\end{figure}

We solve this problem using the $f$-wave solver developed in \cite{leveque2002}.
\Fig{stego_comp} shows a comparison of results using Clawpack and 
SharpClaw on this problem, with 24 cells per layer.
The SharpClaw results are significantly more accurate.

\begin{figure}
\centerline{
\subfigure[SharpClaw.\label{fig:stego_tr_sc}]{\includegraphics[width=3in]{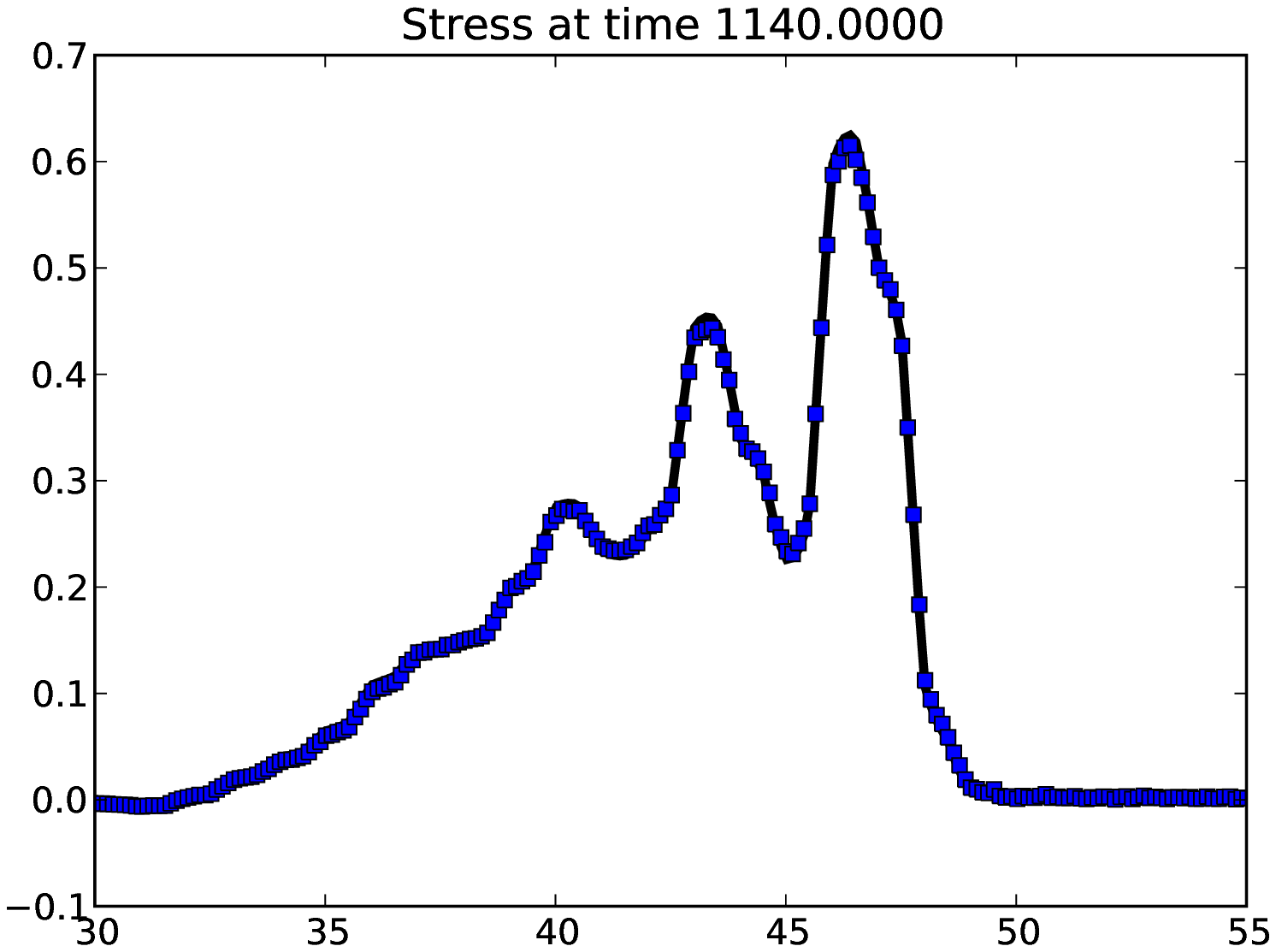}}
\subfigure[Clawpack.\label{fig:stego_tr_cp}]{\includegraphics[width=3in]{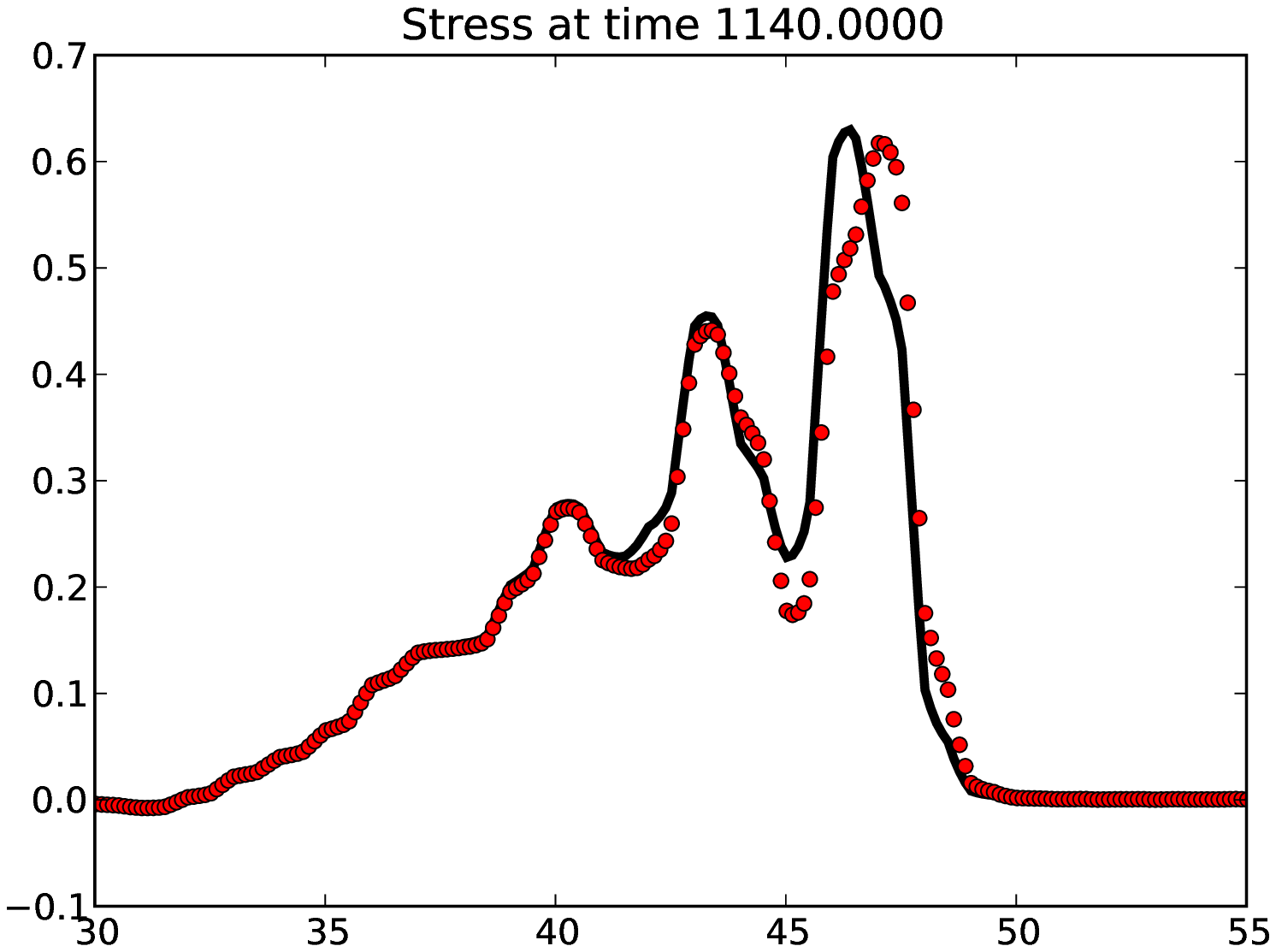}}}
\caption[Comparison of forward solution and time-reversed solution 
stegotons.]{Comparison of forward solution (black line) and 
time-reversed solution (symbols).\label{fig:stego_tr}}
\end{figure}

Solutions of \eqref{nel_pde} are time-reversible in the absence of shocks.  
As discussed in \cite{PhDDKetcheson2009,ketcheson2010}, the effective dispersion
induced by material inhomogeneities suppresses the formation of shocks,
for small amplitude initial and boundary conditions, rendering the
solution time-reversible for very long times.
This provides a useful numerical test.
We solve the stegoton problem numerically up to
time $T$, then negate the velocity and continue solving to time $2T$.  
The solution at any time $2T-t_0$, with $t_0\le T$, 
should be exactly equal to the solution at $t_0$.
We take $T=600$ and $t_0=60$.  \Fig{stego_tr_sc} shows the
solution obtained using SharpClaw on a grid with 24 cells
per layer.  The $t=1140$ solution (squares) is in excellent agreement
with the $t=60$ solution (solid line).  In fact, the maximum point-wise
difference has magnitude less than $2\times10^{-2}$.  Using a
grid twice as fine, with 48 cells per layer, reduces the point-wise
difference to $1\times10^{-3}$.
The Clawpack solution, computed on the same grid (24 cells per
layer), is shown in \Fig{stego_tr_cp}.  Again, the SharpClaw solution is
noticeably more accurate.  For a more detailed study of this time-reversibility test,
we refer to \cite{ketcheson2010}.

%Table \ref{tbl:trconv} lists 
%the maximum pointwise difference (in absolute value) between the solutions
%at $t=60$ and $t=1140$, obtained for a range of grids with both Clawpack 
%and SharpClaw.  By using finer grids,
%the Clawpack solution also appears to converge to the early time solution.
%
%\begin{table} \centering
%\begin{tabular}{r|cc|cc} \hline
%& \multicolumn{2}{c}{Clawpack} & \multicolumn{2}{c}{SharpClaw} \\
%N & $L^\infty$ error & Rate & $L^\infty$ error & Rate \\ \hline
%12 & $4.36 \times 10^{-1}$ & -    &  $5.89 \times 10^{-1}$ & -    \\
%24 & $8.15 \times 10^{-2}$ & 2.42 &  $8.14 \times 10^{-3}$ & 6.18 \\
%48 & $1.54 \times 10^{-2}$ & 2.40 &  $9.81 \times 10^{-4}$ & 3.05 \\
%96 & $3.29 \times 10^{-3}$ & 2.23 &  $3.72 \times 10^{-4}$ & 1.40 \\ 
%192& $7.41 \times 10^{-4}$ & 2.15 &  $7.86 \times 10^{-5}$ & 2.24 \\ 
%\end{tabular}
%\caption{Maximum pointwise error in time-reversal test using Clawpack and
%Sharpclaw.
%The quantity $N$ is the number of computational cells per layer of the
%medium.\label{tbl:trconv}}
%\end{table}

\subsection{Shallow Water Flow}
The next example involves solution of the two-dimensional shallow water equations
\begin{subequations} \label{eqn:shallowWater2DCons}
\begin{align}
\label{eqn:shallowWater2DConsMass}
h_t + (hu)_x + (hv)_y & = 0 \\
\label{eqn:shallowWater2DConsMomx}
(hu)_t + \left(\frac{1}{2}hu^2 + \frac{1}{2}gh^2\right)_x + (huv)_y & = -g h b_x \\
\label{eqn:shallowWater2DConsMomy}
(hv)_t  + (huv)_x + \left(\frac{1}{2}hv^2 + \frac{1}{2}gh^2\right)_y & = -g h b_y,
\end{align}
\end{subequations} 
where $h$, $u$ and $v$ are the depth of the fluid and the velocity components
in $x$ and $y$ directions, respectively.  The function $b(x,y)$ is the bottom
elevation and the constant $g$ represents gravitational acceleration.  

Two test cases are presented: a radially symmetric 
dam-break problem over a flat bottom and a small perturbation 
of a steady state over a hump.  A Roe solver with entropy fix is used in both cases.

%%%%%
\subsubsection{Radial dam-break problem}
\label{subsubsec:1Dradial-dam-break-problem}
This problem consists in computing the flow induced by the instantaneous collapse of an idealized circular dam. It is widely used to benchmark numerical methods designed to simulate interfacial flows and impact problems. 

The domain considered is $[-1.25,1.25]\times[-1.25,1.25]$.  The initial depth is
\begin{align}
h(x,y,t=0) & = \begin{cases}
2 & \mbox{ for } \sqrt{x^2+y^2} \le 1/2 \\
1 & \mbox{ for } \sqrt{x^2+y^2} > 1/2,
\end{cases}
\end{align}
and the initial velocity is zero everywhere. 
This tests the ability of the method to compute the 2D propagation of nonlinear waves and the extent
to which symmetry is preserved in the numerical solution. In the presence of
radial symmetry, system (\ref{eqn:shallowWater2DCons}) can be  recast in the following form:
\begin{subequations} 
\begin{align}
\label{eqn:shallowWater1DRadSymMass}
h_t + (hU)_r & = -\frac{hU}{r} \\
\label{eqn:shallowWater1DRadSymMomr}
(hU)_t + \left(\frac{1}{2}hU^2 + \frac{1}{2}gh^2\right)_r & = -\frac{hU^2}{r},
\end{align}
\label{eqn:shallowWater1DRadSym}
\end{subequations}
where $h$ is still the depth of the fluid, whereas $U$ and $r$ are the radial velocity and the radial position. An important feature of these equations is the presence of a source term, which physically arises from the fact that the fluid is spreading out and it is impossible to have constant depth and constant non-zero radial velocity. 

A first comparison between SharpClaw and Clawpack is performed by solving the 1D system 
(\ref{eqn:shallowWater1DRadSym}) on the interval $0\le r \le 2.5$. 
A wall boundary condition and non-reflecting boundary condition are imposed at the left 
and the right boundaries, respectively. The final time for the analysis is taken to be 
$t=1$. The classical $q$-wave Riemann solver based on the Roe linearization is used to solve the Riemann problem at each interface 
(see for instance \cite{levequefvbook} for details), where the left and the right states are computed by using the characteristic-wise WENO reconstruction. The gravitational acceleration 
is set to $g=1$. A highly resolved solution obtained with Clawpack on a grids with $25,600$ cells is used as a reference solution.

\begin{figure}
\centering
\includegraphics[trim=1cm 1cm 1cm 0cm, clip=true, width=5in] {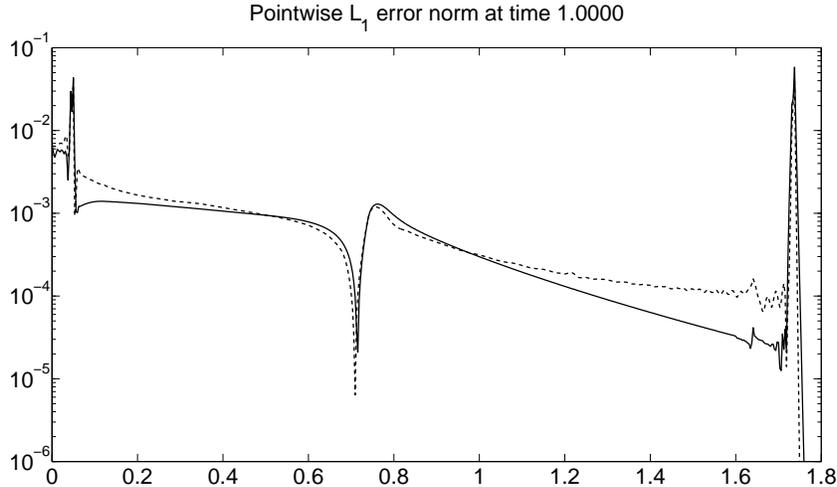}
\caption{Pointwise absolute error for the water height on a grid with $800$ cells. Solid line: SharpClaw
solution; dashed line: Clawpack solution. \label{fig:sw1DRad}}
\end{figure}

\begin{figure}
\centering
\includegraphics[trim=1cm 1cm 1cm 0cm, clip=true, width=5in]{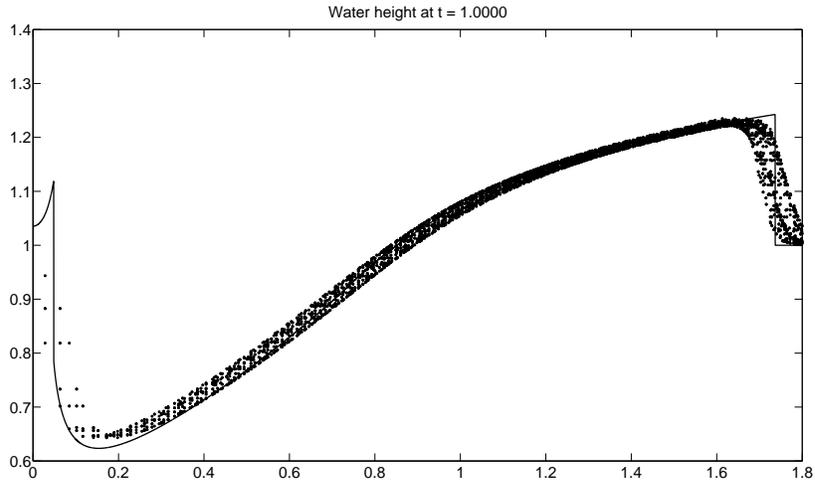}
\caption{Solution for the 2D radial dam-break problem on a grid with $125 \times 125$ cells,
plotted as a function of radius. \label{fig:sw2DRadScatter125}}
\end{figure}

It is well-known that high order convergence is not observed in the presence of shock waves \cite{Majda-1977} and typically the convergence rate is no greater than first order. However, if we plot the difference between the computed solution available at the cell's center and the reference solution conservatively averaged on the same grid, i.e. $E_i = |Q_i - \bar{Q}_i|$, then we can visualize where the errors are largest as well as their spatial structure. Figure \ref{fig:sw1DRad} shows this difference for the water height ($h$) on a grid with $800$ %and $3,200$ 
cells. The reference solution at $t=1$ is shown by the solid line in Figure \ref{fig:sw2DRadScatter125}.  The largest errors in both solutions are near the shocks. In the smooth regions,
the SharpClaw solution is more accurate than that of the Clawpack code.

%\subsubsection{2D Radial Dam-break}
Next we consider the same problem using the full 2D equations 
\eqref{eqn:shallowWater2DCons}.
The SharpClaw and Clawpack codes are tested on two grids with $125 \times 125$ and $500 \times 500$ cells. The final time for the analysis is again taken to be $t=1$.

Figure \ref{fig:sw2DRadScatter125} shows the water height $h$ computed with 
SharpClaw at each cell's center and $t=1.0$ as a function of the radial 
position. The radius is measured respect to the center of the initial 
condition. The 1D reference solution used before is also plotted for comparison. 
It can be seen that the scheme preserves a good radial symmetry, though it 
cannot resolve the shock near the origin. The grid is in fact too coarse. 
Clawpack results are not shown in this figure but indicate similar accuracy and
similarly good symmetry. 

\begin{figure}
\centering
\includegraphics[trim=1.0cm 1cm 1cm 0cm, clip=true, width=5in]{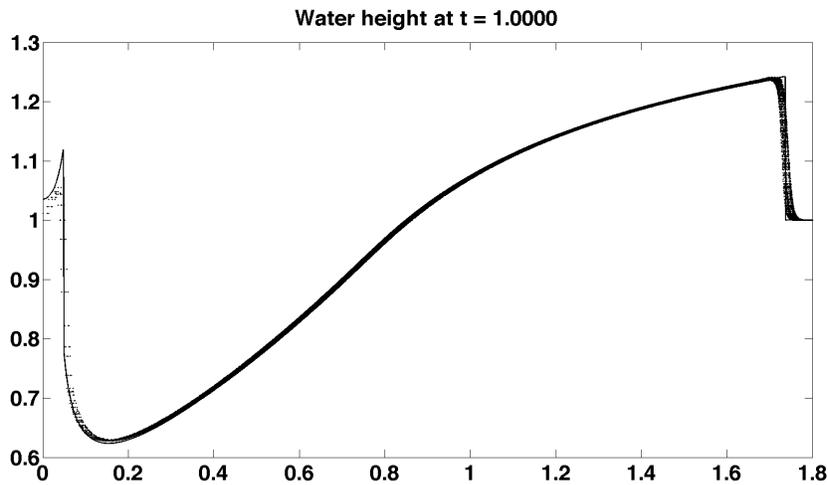}
\caption{Solution for the 2D radial dam-break problem on a grid with $500 \times 500$ cells,
plotted as a function of radius. \label{fig:sw2DRadScatter500}}
\end{figure}

The solutions obtained on the finer grid ($500 \times 500$ cells) are shown in Figure \ref{fig:sw2DRadScatter500}. The effect of the grid refinement is clearly visible. In fact, the solutions gets close to the reference solution. However, the density of the grid near the origin is still too coarse to resolve the shock near the origin to high accuracy.

%%%%
\subsubsection{Perturbation of a steady state solution}
\label{subsubsec:2Dperturb-SWEs}
Conservation laws with source terms often have steady states in which the flux gradient are non-zero but exactly balanced by source terms. A good numerical scheme should be able to preserve such steady states and accurately model small perturbations around these conditions. A classical benchmark test case to investigate these properties is the small perturbation of a 2D steady state water given by LeVeque \cite{LeVeque-1998-QuasiSteady}.

System (\ref{eqn:shallowWater2DCons}) is solved in a rectangular domain $[0,2] \times [0,1]$, with a bottom topography characterized by an ellipsoidal gaussian hump:
$$
b(x,y) = 0.8 \exp(-5(x-0.9)^2 - 50(y-0.5)^2).
$$
The surface is initially flat with $h(x,y,0)=1-b(x,y)$ except for $0.05<x<0.15$, where $h$ is perturbed upward by $\epsilon = 0.01$. The initial discharge in both direction is zero, i.e. $hu(x,y,0)=hv(x,y,0) = 0$. Non-reflecting (i.e., zero-extrapolation) conditions are imposed at all boundaries. The gravitational acceleration is set to $g=9.81$.

\begin{figure}
\centering
\subfigure[Component-wise reconstruction.\label{fig:noBalanced}]{\includegraphics[trim=1.5cm 3cm 1cm 3cm, clip=true, width=3.2in]{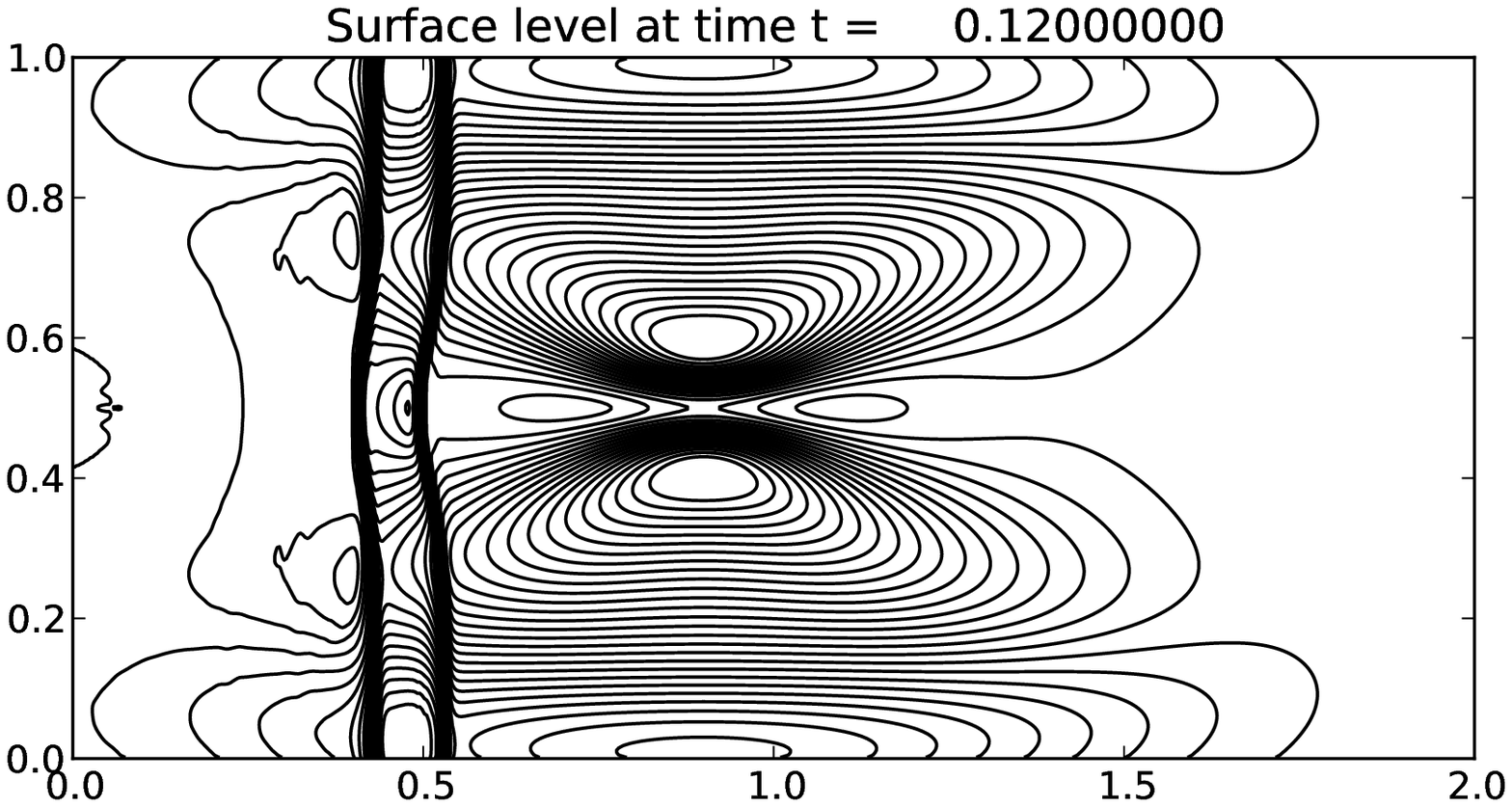}}
\subfigure[$f$-wave-slope reconstruction.\label{fig:balanced}]{\includegraphics[trim=1.5cm 3cm 1cm 3cm, clip=true, width=3.2in]{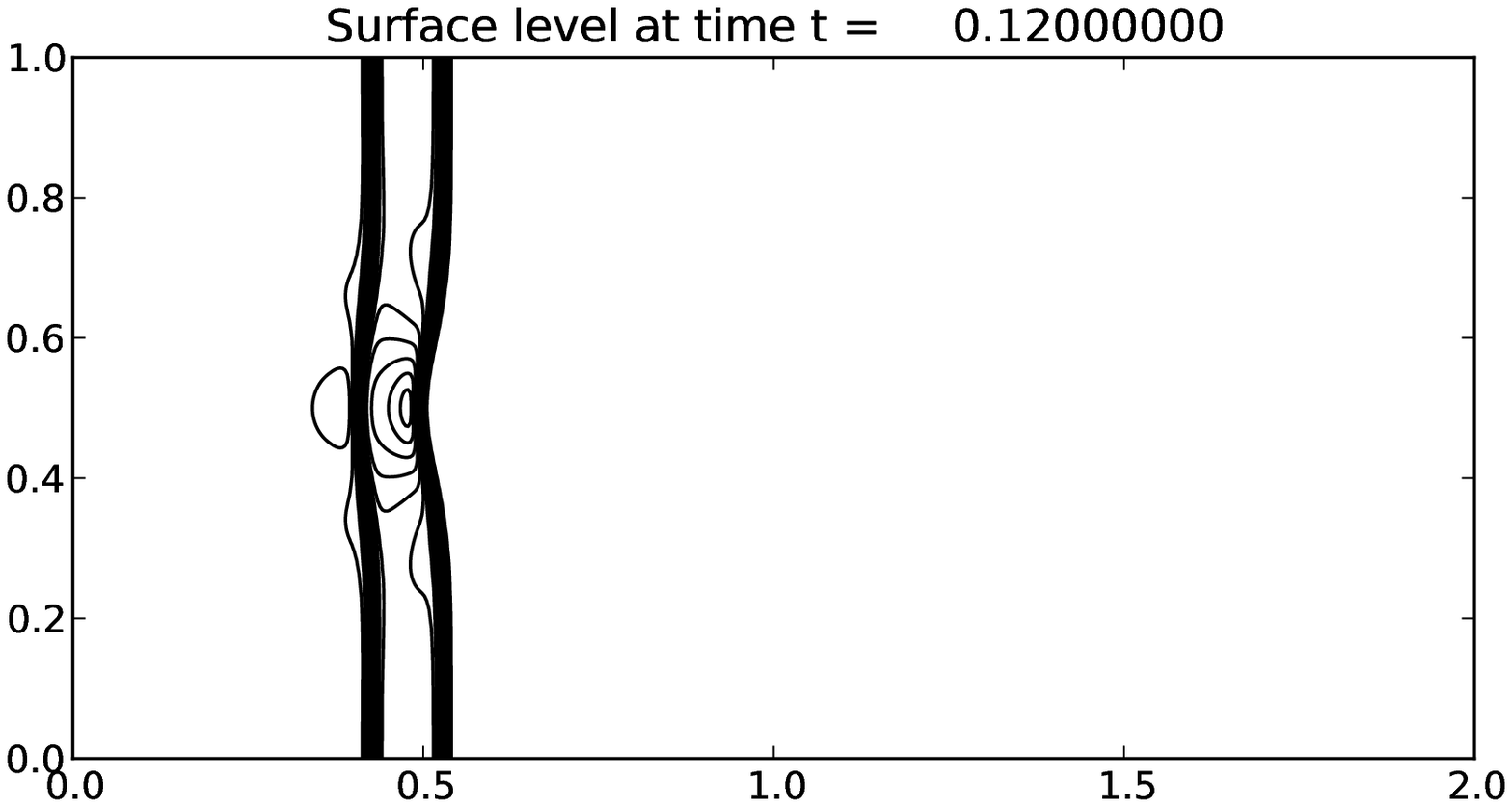}}
\caption{Contour of the surface level $h + b$ at time $t=0.12$ computed with component-wise reconstruction and $f$-wave-slope reconstruction. Contour levels: $0.99942 : 0.000238: 1.00656$.\label{fig:comp-noBalanced-balanced}}
\end{figure}

An effort was made to achieve a well-balanced scheme using the $f$-wave approach
combined with component-wise or characteristic-wise WENO reconstruction, but
this was unsuccessful.  This is not surprising, since the algorithm begins by
reconstructing a non-constant function.
Figure \ref{fig:noBalanced} shows the contour levels of the solution at $t=0.12$ on a fine grid with $600 \times 300$ cells, obtained with the $f$-wave Riemann solver and the component-wise reconstruction approach as a building block for the WENO scheme. The scheme is not well-balanced and spurious waves are generated around the hump. Similar results are obtained using characteristic-wise reconstruction.

\begin{figure}
\centering
\includegraphics[trim=1.0cm 1cm 1cm 0cm, clip=true, width=5in]{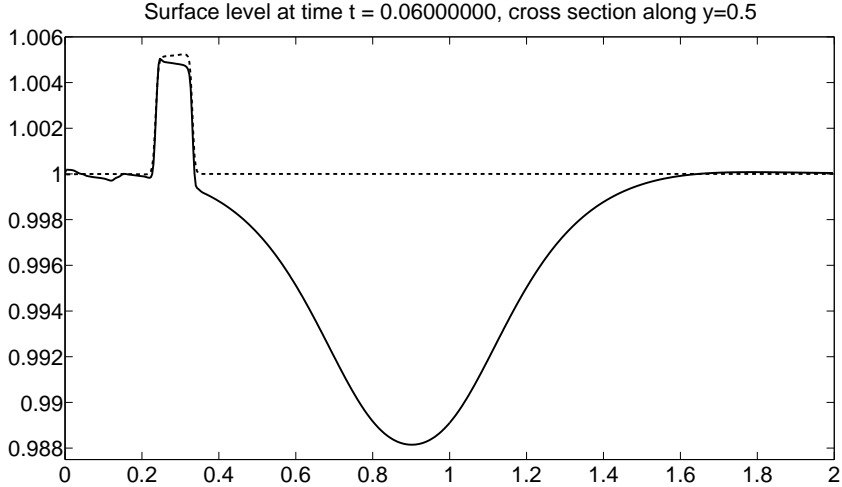}
\caption{Surface level $h + b$ along a cross section at $y=0.5$ and time $t=0.06$. Solid line: $f$-wave Riemann solver and component-wise reconstruction; dashed line: $f$-wave Riemann solver and $f$-wave-slope reconstruction.\label{fig:sw2DSlice-y0p5}}
\end{figure}

In order to balance the scheme, the $f$-wave-slope reconstruction introduced in 
Section \ref{subsec:waveSlopeReconstruction} is used instead. In this approach, the WENO reconstruction is applied to waves computed by solving the Riemann problem
at the cell's interface with the $f$-wave solver.  The bathymetry is approximated by
a piecewise-constant function so that its effect is concentrated at the cell interfaces.
When the source term is included in these Riemann
problems, the resulting waves vanish as shown in Figure \ref{fig:balanced}. In Figure \ref{fig:sw2DSlice-y0p5} the surface level a cross section along $y=0.5$ at time $t=0.06$ computed with both reconstruction approaches (and the $f$-wave Riemann solver) on a uniform mesh with $600 \times 300$ is plotted. This comparison illustrates the different nature of the two approaches. The $f$-wave-slope reconstruction method keeps the
surface flat, whereas the component-wise reconstruction introduces spurious waves which have an amplitude of the order of the disturbance that we want to resolve. 

\begin{figure}
\centering
\subfigure{\includegraphics[trim=1.5cm 3cm 1cm 3cm, clip=true, width=2.3in]{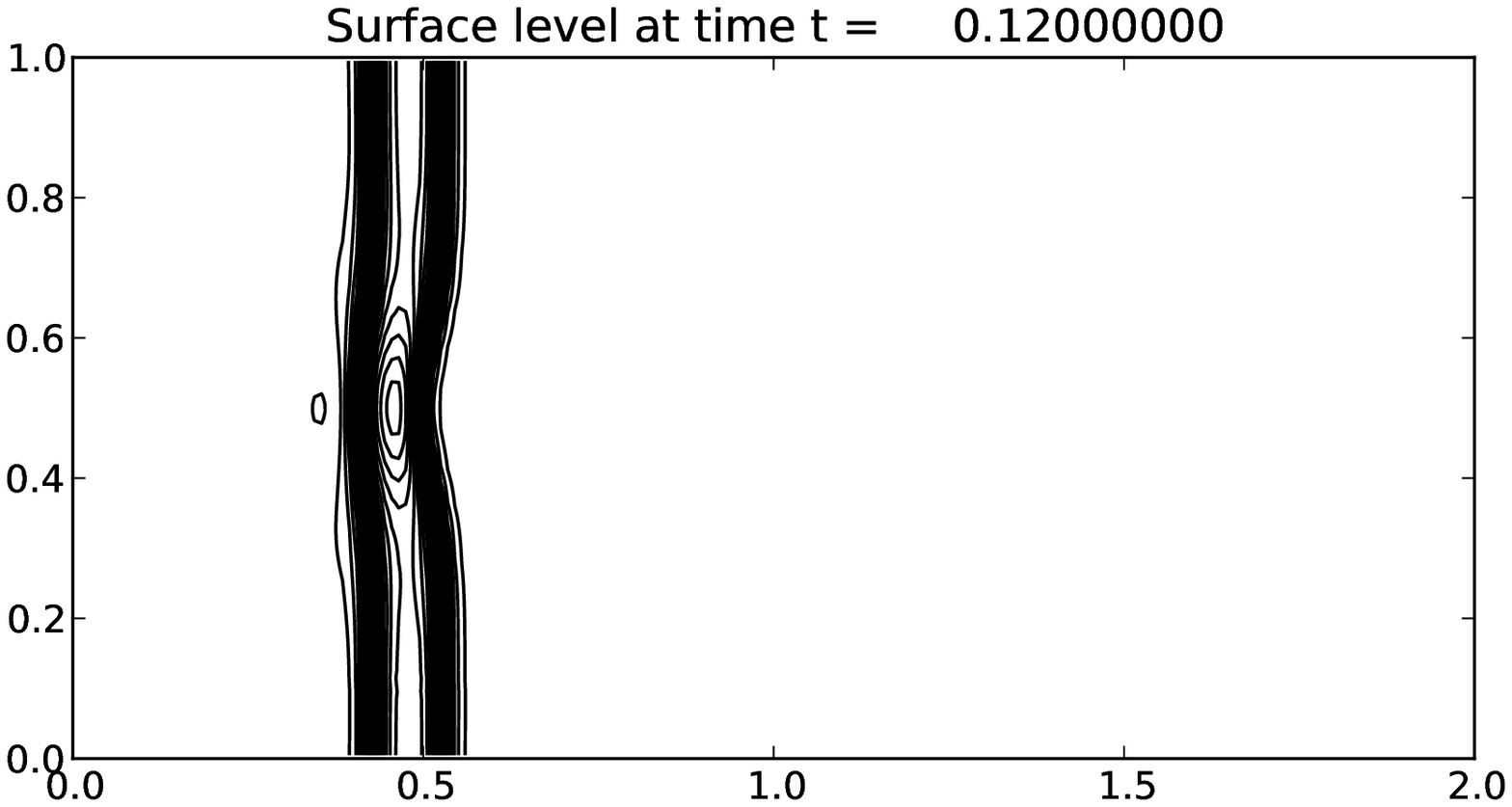}}
\subfigure{\includegraphics[trim=1.5cm 3cm 1cm 3cm, clip=true, width=2.3in]{figures/Pert2DSW-600x300-t012.eps}}
\subfigure{\includegraphics[trim=1.5cm 3cm 1cm 3cm, clip=true, width=2.3in]{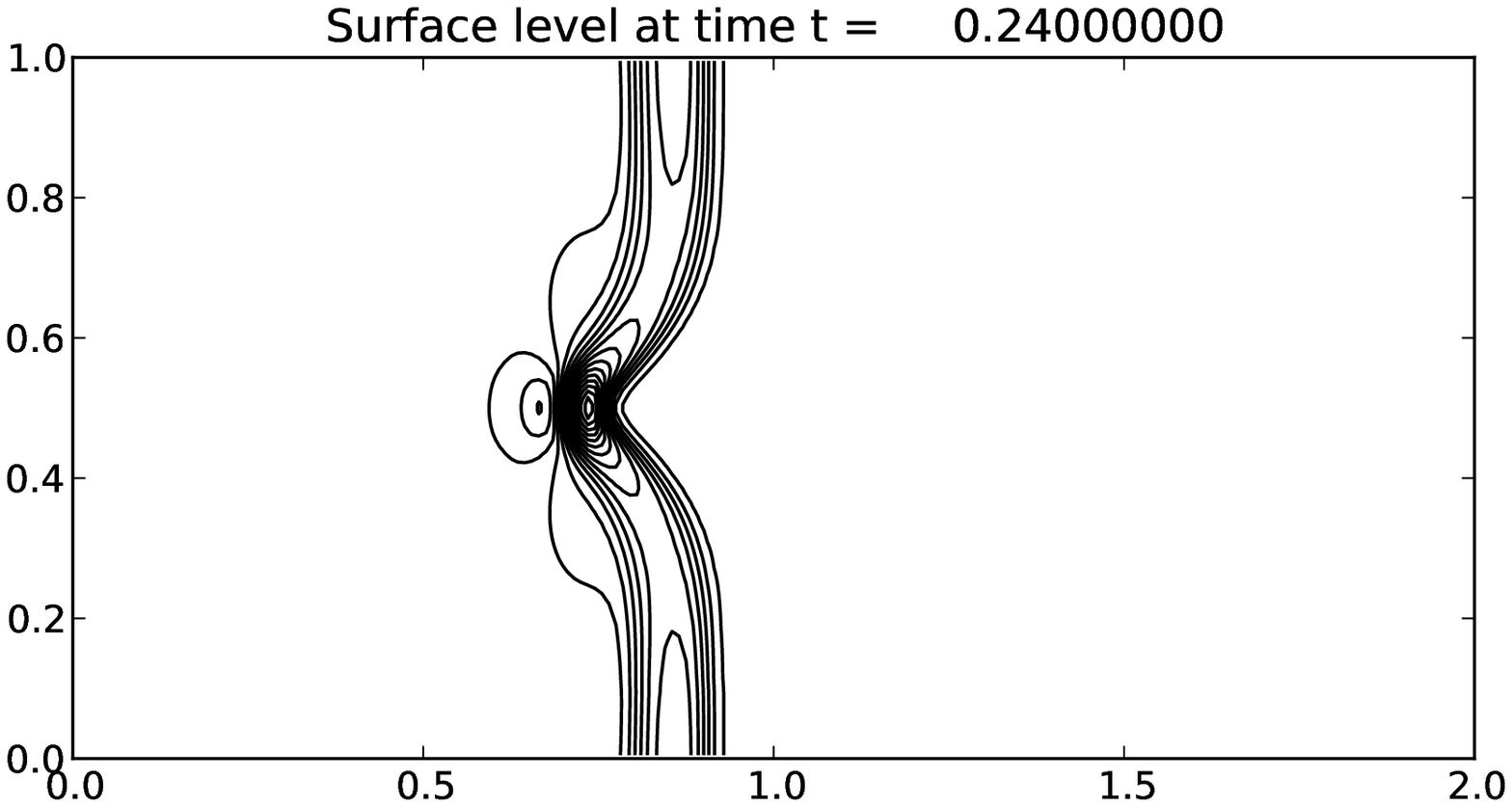}}
\subfigure{\includegraphics[trim=1.5cm 3cm 1cm 3cm, clip=true, width=2.3in]{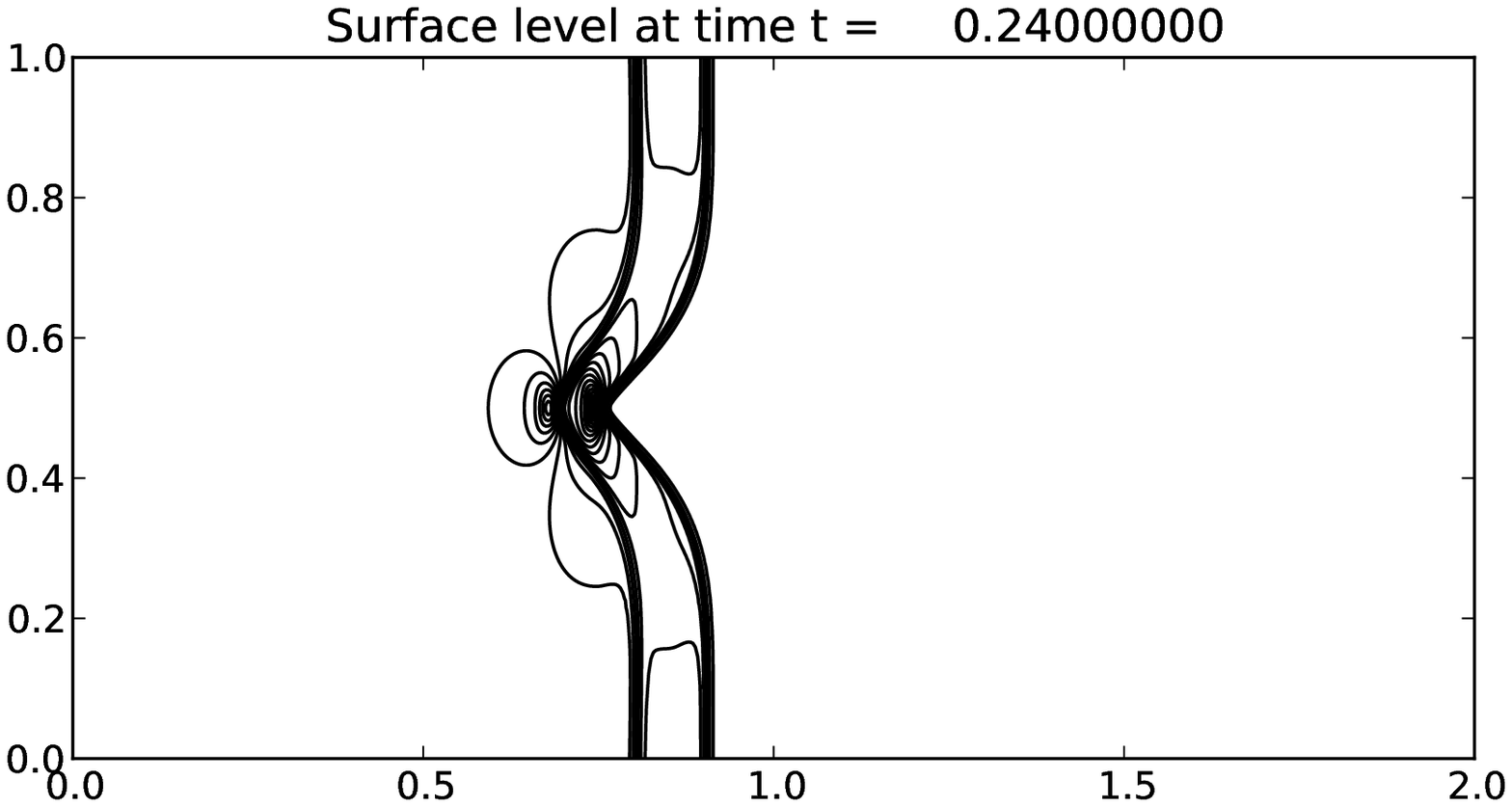}}
\subfigure{\includegraphics[trim=1.5cm 3cm 1cm 3cm, clip=true, width=2.3in]{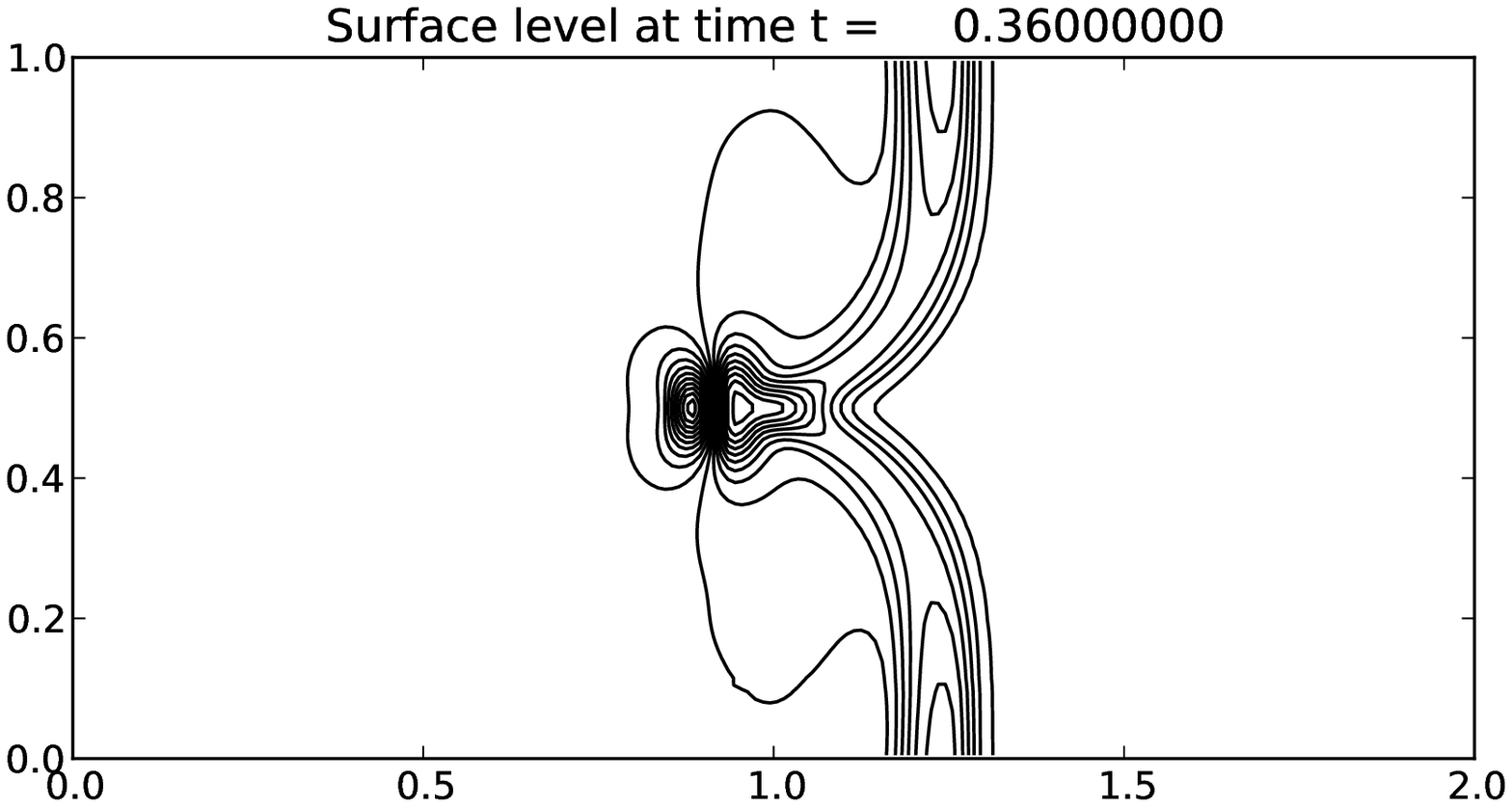}}
\subfigure{\includegraphics[trim=1.5cm 3cm 1cm 3cm, clip=true, width=2.3in]{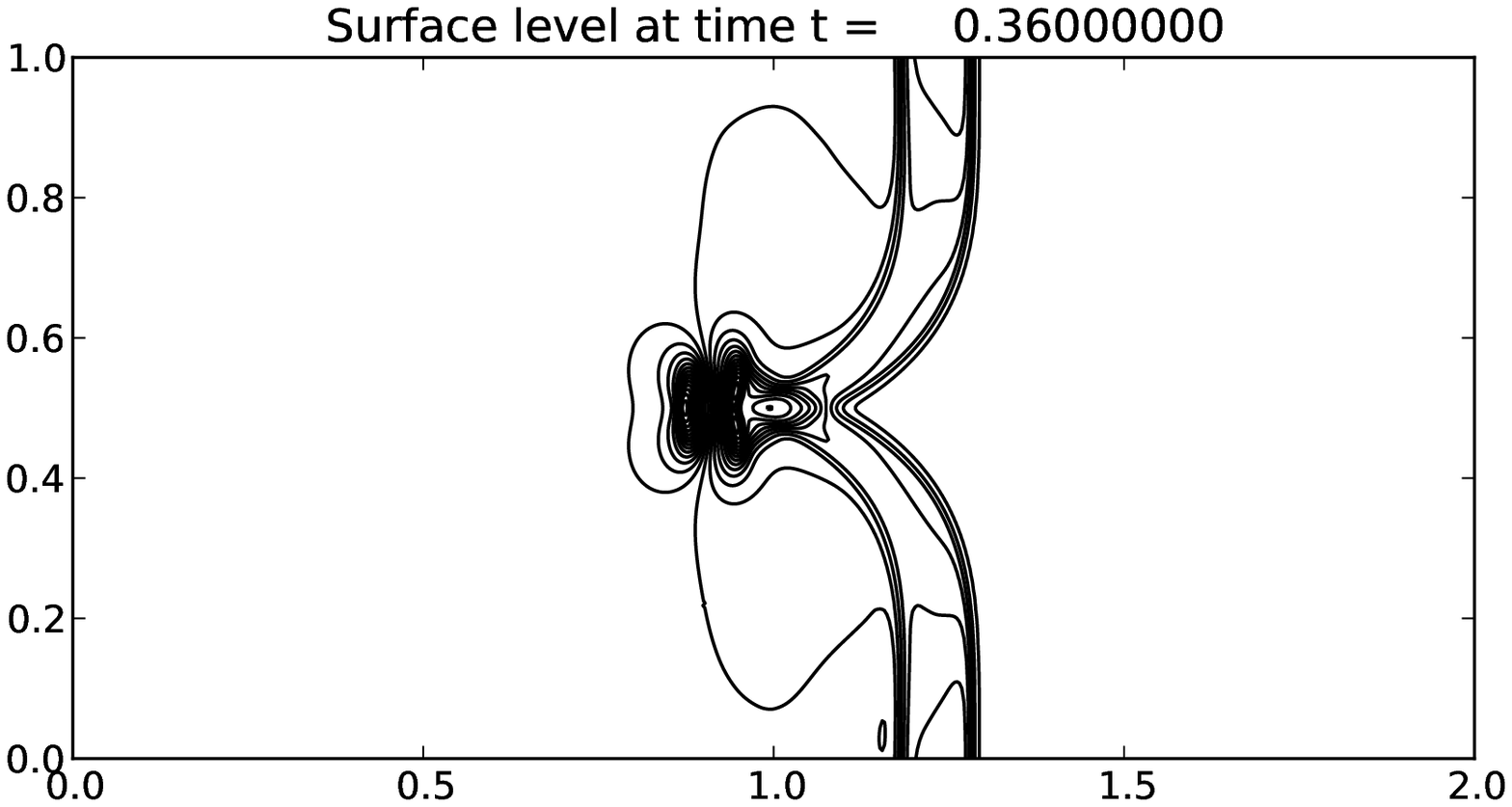}}
\subfigure{\includegraphics[trim=1.5cm 3cm 1cm 3cm, clip=true, width=2.3in]{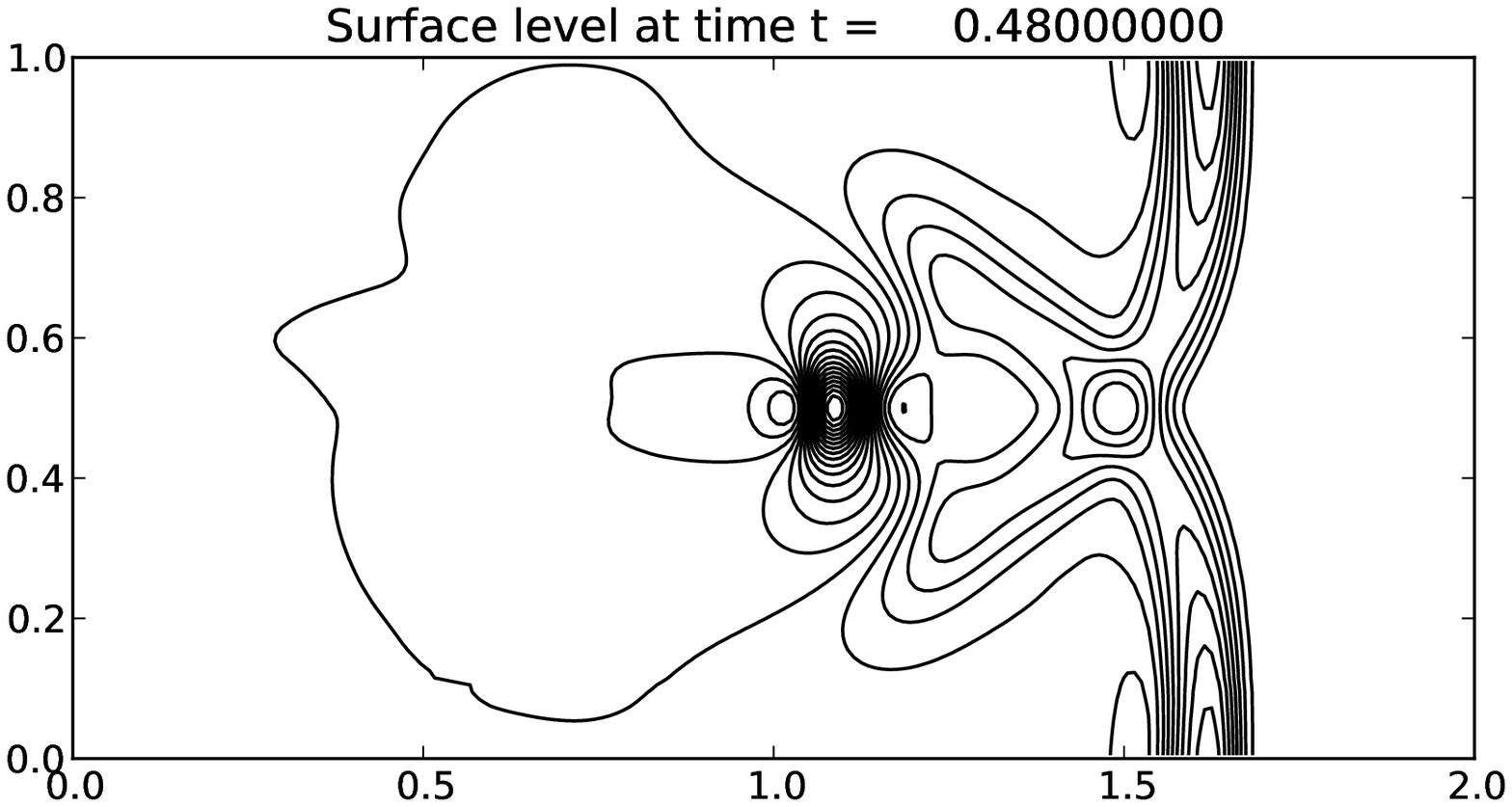}}
\subfigure{\includegraphics[trim=1.5cm 3cm 1cm 3cm, clip=true, width=2.3in]{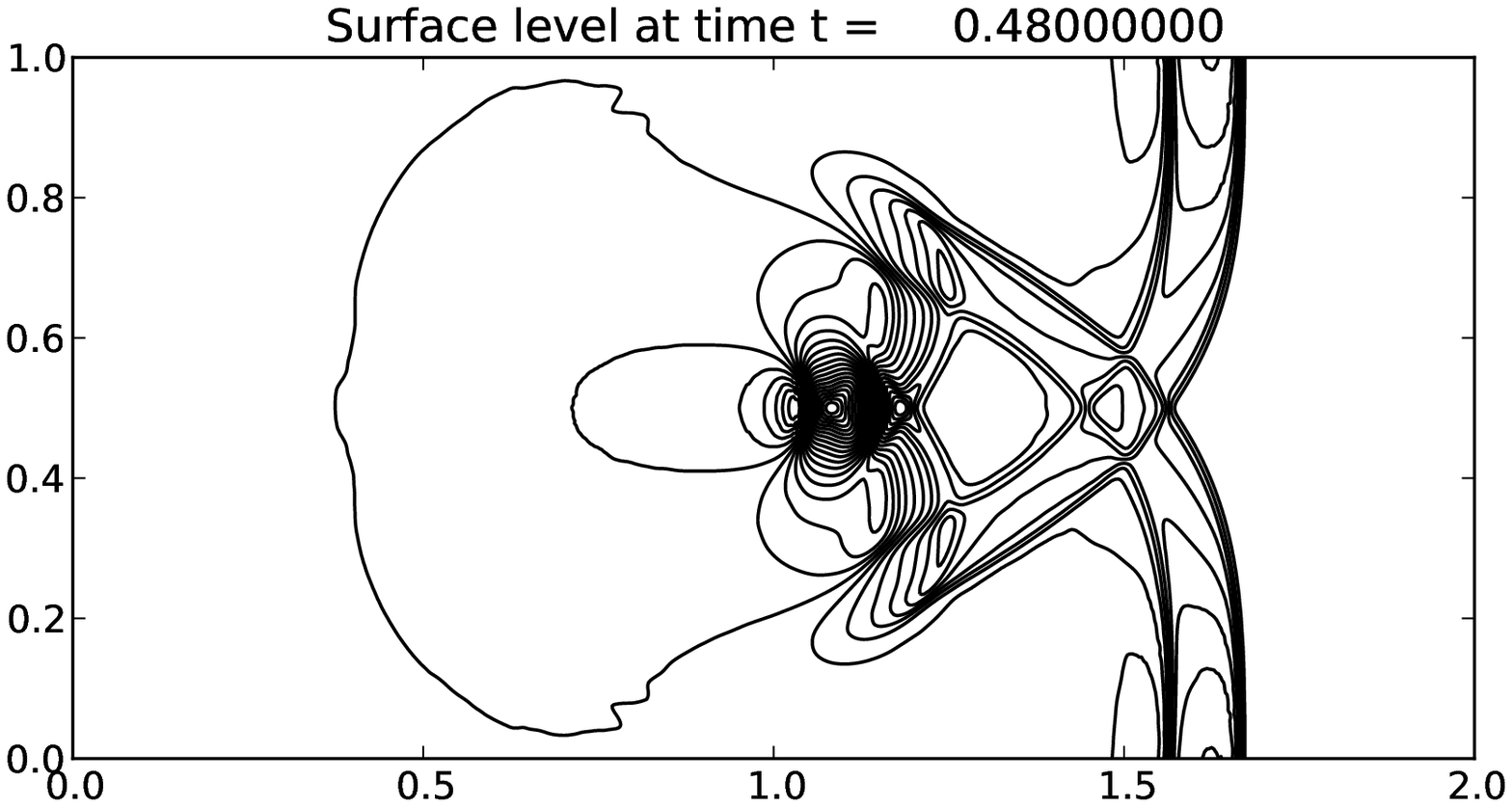}}
\subfigure{\includegraphics[trim=1.5cm 3cm 1cm 3cm, clip=true, width=2.3in]{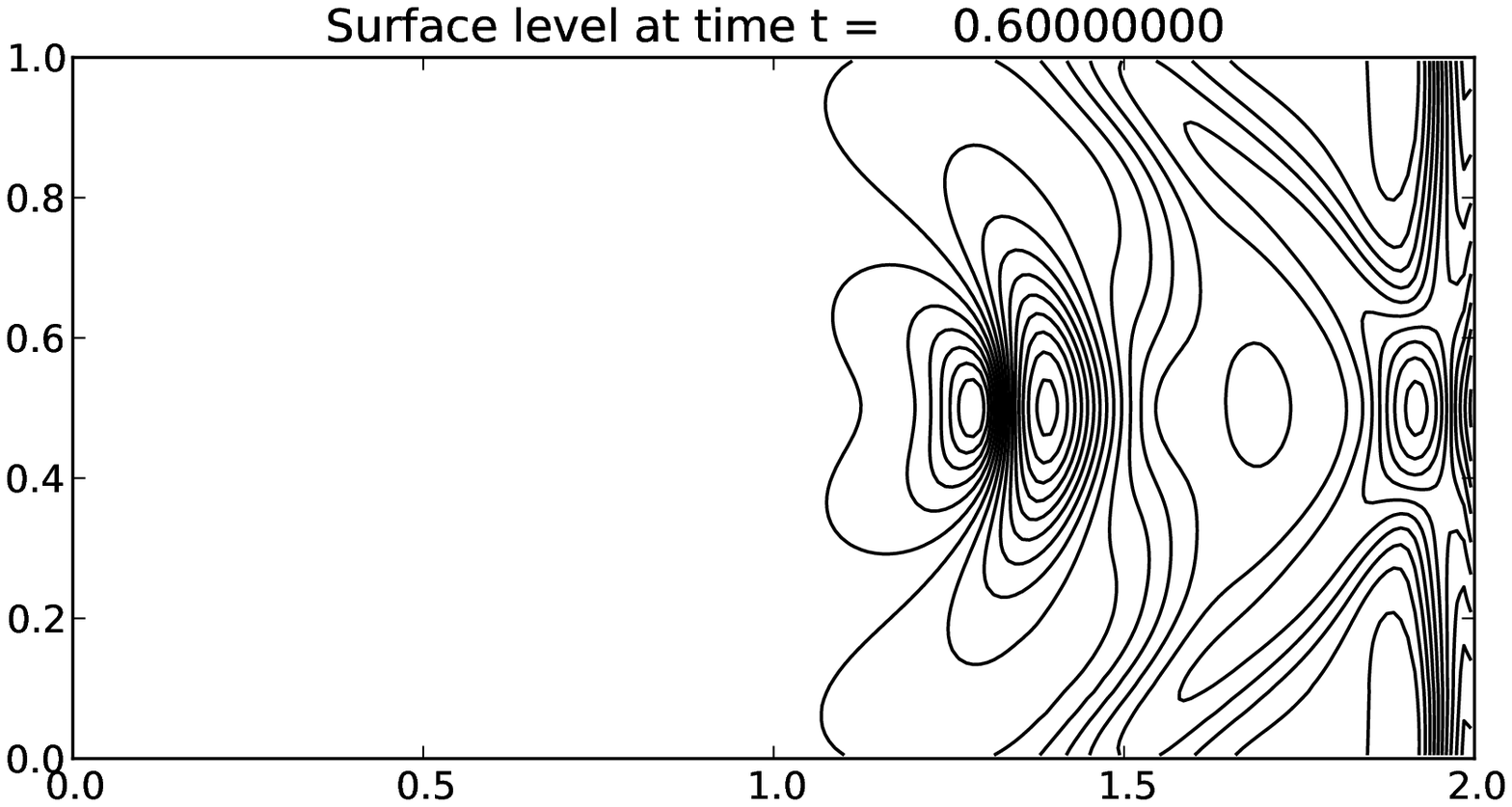}}
\subfigure{\includegraphics[trim=1.5cm 3cm 1cm 3cm, clip=true, width=2.3in]{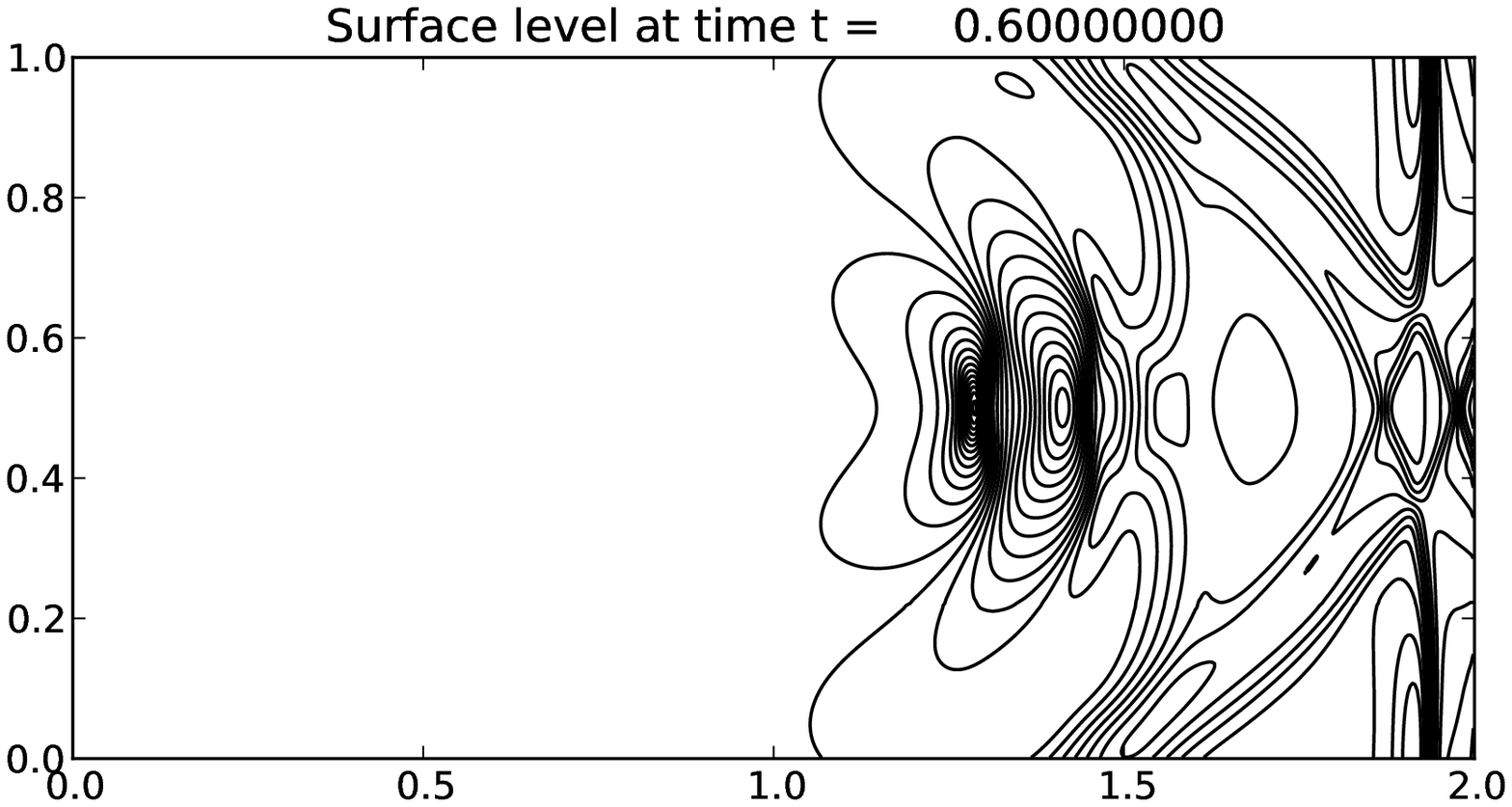}}
\caption{Contour of the surface level $h + b$. $f$-wave-slope reconstruction. 30 uniformly spaced contour lines. $t = 0.12$ from $0.99942$ to $1.00656$; $t = 0.24$ from $0.99318$ to $1.01659$; $t = 0.36$ from $0.98814$ to $1.01161$; $t = 0.48$ from $0.99023$ to $1.00508$; $t = 0.6$ from $0.995144$ to $1.00629$. Left: results with a $200\times100$ cells. Right: results with a $600\times300$ cells. \label{fig:Pert2DSWE-coraseFineGrids}}
\end{figure}

Figure \ref{fig:Pert2DSWE-coraseFineGrids} shows the solution on two uniform meshes with $200 \times 100$ cells and $600 \times 300$ cells, computed using the $f$-wave-slope reconstruction approach. 
The results clearly indicate that the detailed structure of the evolution of such a small perturbation is resolved well even with the relatively coarse mesh. These results agree with those reported in \cite{LeVeque-1998-QuasiSteady}. 
%\captionsetup[subfigure]{labelformat=empty}

In order to investigate the accuracy of our scheme for smooth solutions we also 
have performed a convergence study at a fixed CFL of 0.3 for the same 2D problem.
The smooth initial perturbation is given by 
\begin{align} \label{smooth}
h(x,y,t=0) = \exp(-50(x-0.1)^2)/100
\end{align}
The results of the convergence study are shown in Table \ref{tab:convergence-pert-shallow}.
A highly resolved numerical simulation computed with Clawpack on a $40,000\times 20,000$
grid has been used as a reference solution. 
Since the bathymetry $b(x,y)$ is 
approximated by a piecewise-constant function, both discretizations are formally only
second-order accurate, and the results are roughly consistent with this. 
Nevertheless, the SharpClaw discretization yields significantly more accurate results.

\begin{table}[th]
\caption{Convergence results for the smooth initial condition \eqref{smooth}.  $L_2$-norm of the error as a function of the grid spacing.
\label{tab:convergence-pert-shallow}}
\begin{center}
\begin{tabular}{c|cc|cc}
& \multicolumn{2}{c|}{SharpClaw}& \multicolumn{2}{c}{Clawpack}\\ \hline 
$\Dx=\Dy$ & Error & Order & Error & Order\\ \hline 
 1/10   &  1.14e-2 &      & 2.78e-2 &      \\
\shadeRow  1/20   &  7.00e-3 & 0.70 & 2.09e-2 & 0.41 \\
 1/100  &  8.11e-4 & 1.34 & 6.33e-3 & 0.74 \\
\shadeRow  1/200  &  3.24e-4 & 1.32 & 2.40e-3 & 1.40 \\
 1/1000 &  2.93e-5 & 1.49 & 1.43e-4 & 1.75 \\
\shadeRow  1/2000 &  5.94e-6 & 2.30 & 3.81e-5 & 1.91 \\
\hline
\end{tabular}
\end{center}
\end{table}

\section{Conclusions}
We have presented a general approach to extending the finite volume wave propagation 
algorithm to high order accuracy in one and two dimensions.  
The algorithm is based on a method-of-lines
approach, wherein the semi-discrete scheme relies on high order reconstruction
and computation of fluctuations, including a {\em total fluctuation} term arising
inside each cell.  By using WENO reconstruction and strong stability preserving
time integration, high order accurate non-oscillatory results are obtained,
as demonstrated through a variety of test problems.

This algorithm has several desirable features.  Like the second-order wave propagation
algorithms in Clawpack \cite{leveque1997}, it is applicable to hyperbolic PDEs including linear
nonconservative systems and nonlinear systems with spatially varying flux function.
It has been shown to achieve high order accuracy even for problems with discontinuous
coefficients.  Finally, the algorithm can be adapted to give a well-balanced scheme
for balance laws by use of the $f$-wave approach and a new {\em wave-slope reconstruction} technique.

Hyperbolic systems of equations with both smooth and non-smooth solution have
been used to test the properties and the capabilities of the proposed method.
The schemes have been compared for linear acoustics and nonlinear elasticity
problems in heterogeneous media and for the shallow water
equations with and without bottom topography. Two types of Riemann solver have
been used, i.e the classical ($q$-) wave algorithm and the $f$-wave approach.
The new scheme performed well for all the test cases. 
It gives significantly better accuracy than Clawpack (on the same grid) for smooth problems.

%In two dimensions, although the presented dimension-by-dimension reconstruction WENO scheme is formally second-order accurate for nonlinear systems of equations (see the recent paper of Zhang et al. \cite{Zhang-2011-WENO}),  it is more more accurate than the second-order scheme implemented in Clawpack, at least for the test problems experimented in this work. In addition, if the wave-slope reconstruction algorithm is used as a building block for the WENO scheme, the method is well-balanced and for a small perturbation of the steady state solution of the shallow water equation it can capture very well the small features of the flow. 

A drawback of our implementation of well-balancing is that it requires the
effect of the source term to be approximated entirely at the cell interfaces.
For shallow water equations with smooth bathymetry, this will reduce the 
formal accuracy to second order.  Future work might explore the implementation
of higher order accurate well-balancing for polynomial source terms using 
high order quadrature.
In any case, in two dimensions, the presented dimension-by-dimension reconstruction 
approach is formally only second-order accurate
(see \cite{Zhang-2011-WENO}); however, it gives improved accuracy over the second-order 
scheme implemented in Clawpack for the test problems considered.
Further investigation of different approaches to multidimensional reconstruction 
for problems containing both shocks and rich smooth flow structures is a topic 
of future research.

%Although the work of Zhang et al. \cite{Zhang-2011-WENO} has gave some guidance in the application of high order finite
%volume schemes for simulating shocked flows, further investigation of the accuracy of WENO schemes  In addition, both structured and unstructured grid based high order methods will co-exist and excel for different type of flow problems, a thorough comparison of the WENO scheme with other emerging high order compact schemes for unstructured mesh such as discontinuous Galerkin and spectral difference method will be perform in future work. 
\subsection*{Acknowledgments} The authors are grateful to anonymous referees whose suggestions
improved this paper.  This work was supported in part by NSF grants
DMS-0609661 and DMS-0914942.
\bibliography{sharpclaw}

\end{document}